\documentclass[10pt]{amsart}

\usepackage{amsfonts,amssymb}
\newtheorem{de}{Definition}[section]
\newtheorem{lm}[de]{Lemma}
\newtheorem{pr}[de]{Proposition}
\newtheorem{co}[de]{Corollary}
\newtheorem{re}[de]{Remark}

\newtheorem{te}[de]{Theorem}
\newtheorem{ex}[de]{Example}

\input xy
\xyoption{all}
\begin{document}
\title[Duals of pointed Hopf algebras]{Duals of Pointed Hopf Algebras}
\author{M. Beattie}
\address{Department of Mathematics and Computer Science\\
Mount Allison University
\\ Sackville, N.B., Canada E4L 1E6\\ email: mbeattie@mta.ca}
\thanks{Research partially supported by NSERC}

\date{}
\maketitle

\begin{abstract}

In this paper, we  study the duals of some finite dimensional
pointed Hopf algebras  working over an algebraically closed field
$k$ of characteristic 0. In particular, we study pointed Hopf
algebras with coradical $k[\Gamma]$ for $\Gamma$ a finite abelian
group, and with  associated graded Hopf algebra of the form
${\mathcal B}(V) \# k[\Gamma]$ where  ${\mathcal B}(V)$ is the
Nichols algebra of
 $V= \oplus_i V^{\chi_i}_{g_i} \in
{^{k[\Gamma]}_{k[\Gamma]}}{\mathcal YD}$.   As a corollary to a
general theorem on duals of coradically graded Hopf algebras, we
have that  the dual of ${\mathcal B}(V) \# k[\Gamma]$
  is ${\mathcal
B}(W)\#k[\hat{\Gamma}]$ where $W= \oplus_i W^{g_i}_{\chi_i} \in
^{k[\hat{\Gamma}]}_{k[\hat{\Gamma}]}{\mathcal YD}$. This description
of the dual is used to explicitly describe the Drinfel'd double of
${\mathcal B}(V) \#k[\Gamma]$. We also show that the dual of a
nontrivial lifting $A$ of ${\mathcal B}(V)\#k [\Gamma]$ which is not
itself a Radford biproduct,  is never pointed.
 For $V$ a quantum linear
space of dimension 1 or 2, we describe the duals of some liftings of
${\mathcal B}(V) \# k[\Gamma]$.  We conclude with some examples
where we determine all the irreducible finite-dimensional
representations of a lifting of ${\mathcal B}(V) \# k[\Gamma]$ by
computing the matrix coalgebras in the coradical of the dual.
\end{abstract}

\section{Introduction and Preliminaries}\label{intro}

Unless otherwise specified, we work over $k$, an algebraically closed field of
characteristic 0. The object of this paper is to find the
structure of dual Hopf algebras to some finite dimensional
pointed Hopf algebras with coradical $k[\Gamma]$, with $\Gamma$ a
finite abelian group. A method for classifying pointed Hopf
algebras has been proposed by N. Andruskiewitsch and H.-J.
Schneider; an up-to-date survey of this project may be found in
\cite{assurvey}. Roughly, this approach to studying pointed Hopf
algebras is based on the following observations.

Let $A$ be a Hopf algebra such that the coradical $A_0$ is a Hopf
subalgebra of $A$ so that the coradical filtration of $A$ is a Hopf
algebra filtration. Let $H = gr A$ be the associated (coradically)
graded Hopf algebra \cite[Chapter XI]{Sweedler}, i.e., $H =
\oplus_{n \geq 0} H(n)$ with $H(0) = A_0$ and $H(n) = A_n / A_{n-1}$
for $n \geq 1$. (Recall from \cite{cm} that a graded  coalgebra $C=
\oplus^\infty_{i=0} C(i)$ is coradically graded if $C(0) =C_0$ and
$C_1 =C(0) \oplus C(1)$, where $C_m$ denotes the coradical
filtration and also if $C_0$ is one-dimensional then $C$ is called
strictly graded.) Then $\pi$, the Hopf algebra projection of $H$
onto $H_0 = A_0$ \cite[5.4.2]{mont} splits the inclusion $i: H_0 \to
H$ and thus $H \cong {\mathcal R} \# H_0$, the Radford biproduct of
${\mathcal R}$ and $H_0$,
 where ${\mathcal R} = H^{co \pi} = \{x \in H : (id \otimes \pi )\triangle(x)
 = x \otimes 1\}$ is the set of
$\pi$-coinvariants.  $\mathcal R$ is a Hopf algebra in
$^{H_0}_{H_0}{\mathcal YD}$, the category of Yetter-Drinfel'd
modules over $H_0$. In particular, if $H_0 = k[\Gamma]$, for
$\Gamma$ a finite abelian group, as will be the case in the examples
of interest in this paper, then $^{k[\Gamma]}_{k[\Gamma]}{\mathcal
YD}$ is the set of left $k [\Gamma]$-modules with a $\Gamma$-grading
such that the $\Gamma$-action preserves the grading.

 Also $\mathcal R$
inherits a grading from $H$, i.e. $\mathcal R$ is graded by
${\mathcal R}(n) = {\mathcal R} \cap H(n)$ with ${\mathcal R}(0) = k
\cdot 1$, and ${\mathcal R}(1) = P({\mathcal R})$, the space of
primitive elements of $\mathcal R$. If $V = {\mathcal R}(1)$
generates $\mathcal R$ as a Hopf algebra, then ${\mathcal R} =
{\mathcal B}(V)$, the Nichols algebra of $V$. In any case ${\mathcal
B}(V) \subseteq \mathcal R$. Much more detail about Nichols algebras
can be found in \cite{grana} or in the original source material
\cite{nichols}. An excellent comprehensive survey of this approach
with all the definitions and a compendium of results from various
papers in the most generality may be found in \cite{assurvey}.

Now, to construct a finite dimensional pointed Hopf algebra with
coradical $k[\Gamma]$, the idea is to first find
 $V \in ^{k[\Gamma]}_{k[\Gamma]}{\mathcal YD}$ with ${\mathcal
B}(V)$  finite dimensional, and form the coradically graded Hopf
algebra $H = {\mathcal B}(V) \# k[\Gamma]$. This construction is
called bosonization or the Radford biproduct. A Hopf algebra $A$
such that $gr A \cong H$ is called a lifting of $H$ and if $A \not
\cong gr A$, then $A$ is a nontrivial lifting of $H$. If $V$ is a
quantum linear space, or is of Cartan type $A_n$ or of Cartan type
$B_2$, the liftings of ${\mathcal B}(V) \# k[\Gamma]$ have been
computed in \cite{asp3}, \cite{bdg}, \cite{asa2},
\cite{assurvey},\cite{ad}, \cite{bdr}, respectively. We will study
duals of liftings of some quantum linear spaces.

We abbreviate $^{k[\Gamma]}_{k[\Gamma]}{\mathcal YD}$ to
$^{\Gamma}_{\Gamma}{\mathcal YD}$ for the category of
Yetter-Drinfel'd modules over $k[\Gamma]$. For $V \in
^{\Gamma}_{\Gamma}{\mathcal YD}$, we write $V^{\chi}_g$ with $g \in
\Gamma, \chi \in \hat{\Gamma}$ for the set of $v \in V$ with the
action of $\Gamma$ on $v$ given by $h \to v = \chi (h) v$ and the
coaction by $\delta (v) = g \otimes v$. If $\Gamma$ is an abelian
group, and $V \in ^{\Gamma}_{\Gamma}{\mathcal YD}$ then $V =
\oplus_{\scriptstyle g \in \Gamma \atop\scriptstyle \chi \in
\hat{\Gamma}} V^\chi_g$ \cite[Section 2]{assurvey}. If $V=
\oplus^t_{i=1} kv_i$ with $v_i \in V^{\chi_i}_{g_i} $ then  the
braid matrix $B$ of $V$ is $(b_{ij}) =( \chi_j (g_i)).$

Under our assumptions on $k$ and assuming $\Gamma$ is finite
abelian, we identify $k[\Gamma]^*$ with the Hopf algebra
$k[\hat{\Gamma}]$ where $\hat{\Gamma}$ is the character group of
$\Gamma$. Let $\varphi: k[\Gamma]^* \to k[\hat{\Gamma}]$ be the
isomorphism between these Hopf algebras. For $g \in \Gamma$, the
element $g^{**} \in Alg(k[\hat{\Gamma}], k)$ is defined by
$g^{**} (\chi) = \chi(g)$. We usually simply write $g$ for a
character of $\hat{\Gamma}$ instead of $g^{**}$.

Unless otherwise noted, $\Gamma$ will denote a finite abelian
group.

For $A$ a Hopf algebra with $g,h \in G(A)$, the group of
grouplike elements of $A$, then $P(A)_{g,h} = \{ x \in A |
\triangle x = x \otimes g +h \otimes x\}$, the $(g,h)$-primitives
of $A$. For each $g,h$, we write $P'(A)_{g,h}$ for a subspace of
$P(A)_{g,h}$ such that $P(A)_{g,h} =k (g-h) \oplus P'(A)_{g,h}$
\cite[\S5.4]{mont}.

We use the Sweedler summation notation throughout.  $\overline{S}
$ will denote the composition inverse to the antipode $S$.

Throughout, we write $I_j$ for the $ j \times j$ identity matrix.
Also we use the notation ${\mathcal M}^c(r,k)$ for an $r \times r$
matrix coalgebra over $k$ and we call a basis $e(i,j), 1 \leq i,j
\leq r$, a matrix coalgebra basis if $\triangle (e(i,j)) =
\sum^r_{k=1} e(i,k) \otimes e(k,j)$, and $\epsilon (e(i,j)) =
\delta_{ij}$.

Our specific constructions and examples will focus on the case of
$V$ a quantum linear space.

\begin{de}\label{deqls} $V = \oplus^t_{i=1} kv_i \in {^{\Gamma}_{\Gamma}{\mathcal YD}}$
with $v_i \in V^{\chi_i}_{g_i}$ is called a quantum linear space if
$\chi_i (g_j) \chi_j (g_i) =1$ for $i \neq j$, i.e. in the braid matrix
$B$, $b_{ij} b_{ji} =1$ for $i \neq j$.
\end{de}

Recall that for  $V$ a quantum linear space as above, ${\mathcal
B}(V) \# k [\Gamma]$ is the Hopf algebra generated by the $(1,
g_i)$-primitives $v_i \#1$ (usually written just $v_i$) and the
grouplike elements $h \in \Gamma$. Multiplication is given by $hv_i
= \chi_i (h) v_i h$ and $v_i v_j = \chi_j (g_i) v_j v_i$. Also  if
$\chi_i(g_i)$ is a primitive $r_i$th root of unity, $v_i^{r_i} =0$,
and then dim ${\mathcal B}(V) = \Pi^t_{i=1} r_i$.

\begin{pr}\label{liftqls}(see \cite{asp3}
or \cite{bdg}) For $V$ a quantum linear space with $\chi_i(g_i)$
 a primitive $r_i$th root of unity, all liftings $A$ of ${\mathcal
B}(V) \# k [\Gamma]$ are Hopf algebras generated by the
grouplikes and by $(1, g_i)$-primitives $x_i, 1 \leq i \leq t$
where
$$hx_i = \chi_i(h) x_ih;$$
$$x_i ^{r_i} = \alpha_{ii} (g_i^{r_i} -1);$$
$$x_i x_j = \chi_j (g_i) x_j x_i + \alpha_{ij} (g_i g_j -1).$$

We may assume $\alpha_{ii} \in \{ 0,1\}$ and then we must have
that \begin{equation}\label{alphaii} \alpha_{ii} =0 \mbox{ if }
g_i^{r_i} =1 \mbox { or }\chi_i^{r_i} \neq \epsilon \mbox{ and }
\alpha_{ij} =0 \mbox{ if }g_i g_j =1 \mbox{ or }\chi_i \chi_j \neq
\epsilon.
\end{equation}
   Note that $\alpha_{ji} = - \chi_j (g_i)^{-1}
\alpha_{ij} = - \chi_i (g_j) \alpha_{ij}$. Thus the lifting $A$ is
described by a matrix ${\mathcal A} = ( \alpha_{ij})$ with 0's or
1's down the diagonal and with $\alpha_{ji} =- \chi_i (g_j)
\alpha_{ij}$ for $i \neq j$. We associate ${\mathcal A} = (0)$ with
${\mathcal B}(V) \# k [\Gamma]$. \qed
\end{pr}

\section{The dual of a Radford biproduct.}

The first theorem in this section is in the form suggested by N.
Andruskiewitsch who points out that the ideas for this proof are
already in the literature although the statement does not appear
explicitly. (See \cite{asadvances}, \cite{ag}, \cite{nichols}.)
The theorem in the original version of this paper corresponded to
Corollary \ref{cordual} and follows directly from this more
general formulation.
  In fact it is shown in \cite[3.2.30]{ag} that for $V
 \in {^H_H\mathcal{YD}}$ with ${\mathcal B}(V)$ finite dimensional, then
 ${\mathcal B}(^*V) \cong (^* {\mathcal B}(V))^{opcop}$ in
 the category $^H_H\mathcal{YD}$.  Note that the categories
 $^H_H\mathcal{YD}$ and $^{H^*}_{H^*}\mathcal{YD}$ are equivalent for $H$
 finite dimensional with bijective antipode by \cite[2.2.1]{ag}.
 The results of Corollary \ref{cordual} give a special case.

For $k$ any field, let $H$ be a graded Hopf algebra, $H = \oplus_{n
\ge 0} H^n$ with finite dimensional components $H^n$. Let $L =
\oplus_{n \ge 0} (H^n)^*$, with the ``dual'' Hopf algebra structure.
Assume that $H_0 = H^0$, i.e., $H^0$ is cosemisimple. (Note that
without restriction on $k$ this does not imply that $L^0$ is
cosemisimple.) Let $\pi: H \to H^0$, $p: H \to H^1$, $\sigma: L \to
L^0$ be the homogeneous projections. Then, as discussed in Section
\ref{intro}, $H \simeq {\mathcal R} \# H^0$, and  $L \simeq
{\mathcal S} \# L^0$ where ${\mathcal R} = H^{co\pi}$ and ${\mathcal
S} = L^{co \sigma}$.

\begin{lm} For $k$ any field and $H,L$ as above, the following are equivalent:

(i) $H$ is coradically graded.

(ii) $L$ is generated as an algebra by  $L^0$ and $L^1$.

\end{lm}
\begin{proof} Suppose (i) holds. By \cite[2.2]{cm}, $H$ is coradically
graded if and only if $H_1 = H^0 + H^1$.  Now consider the map
$\phi$ from $H^2$ to $H^1  \otimes H^1$ defined by $\phi := p
\otimes p \circ \triangle$. If $\phi(a)=0$, then $\triangle(a) \in
H^0 \otimes H^2 + H^2 \otimes H^0$, and then since $H^0 = H_0$, $a
\in H_1$, which is impossible unless $a=0$ since $H_1 \cap H^2 = 0$.
Thus $\phi$ is injective. By induction (or see the proof of
\cite[Lemma 5.5]{asadvances}), the map $(p \otimes p \otimes \ldots
\otimes p)\circ \triangle^{n-1}$ from $H^n$, the $n$th homogeneous
component of $H$, to $H^1 \otimes \ldots \otimes H^1$ is injective
for $n \geq 2$. Clearly this is equivalent to the surjectivity of
the multiplication map from $L^1 \otimes L^1 \otimes \ldots \otimes
L^1$ to $L^n$.

Conversely, if the multiplication map from $L^1 \otimes L^1 \otimes
\ldots \otimes L^1$ to $L^n$ is surjective, so that the map $(p
\otimes p \otimes \ldots \otimes p)\circ \triangle^{n-1}$ from $H^n$
to $H^1 \otimes \ldots \otimes H^1$ is injective, then since $(p
\otimes p \otimes \ldots \otimes p)\circ \triangle^{n-1}$ acting on
$H_1$ is zero for any $n >1$, then $H_1 \cap H^n = 0$ for all $n
\geq 2$. Thus,  the containment $H^1 + H^0 \subseteq H_1$ is an
equality, and $H$ is coradically graded. \end{proof}

 From the lemma above, we immediately obtain:
\begin{te}\label{dual}Let $k$ be any field and let  $H$ and $L$ be as above,
with the additional assumption that $H^0=H_0$ is cosemisimple and
also semisimple, so that $L_0$ is also. The following are
equivalent:

(i) $H$ is coradically graded and is generated as an algebra by
$H^0$ and $H^1$.

(ii) $L$ is coradically graded and is generated as an algebra by
$L^0$ and $L^1$.

\end{te}

A Hopf algebra $H$ satisfying the conditions of Theorem \ref{dual}
is called a ``bialgebra of type one'' in the original terminology of
Nichols' \cite{nichols}. Then, as briefly outlined in Section
\ref{intro} and explained in detail in \cite[Lemmas 2.4, 2.5]{asp3}
or in the survey \cite[Section 2]{assurvey}, $H \simeq {\mathcal
B}(V) \# H^0$, where $V = {\mathcal R}(1) = P({\mathcal R})$, the
space of primitives of ${\mathcal R} = {\mathcal B}(V)$, and
  $L \simeq {\mathcal B}(W) \# L^0$, where $W = {\mathcal S}(1) = P({\mathcal S})$.
   Here ${\mathcal R}={\mathcal
  B}(V)$ and ${\mathcal S}={\mathcal B}(W)$ are Nichols algebras; in particular, they are strictly graded
  coradically graded Hopf
  algebras in the categories
  $_{H_0}^{H_0}{\mathcal YD}$ and $_{L_0}^{L_0}{\mathcal YD}$
  respectively and are generated as
  algebras by their 1-graded component.

There is a relation between  the Yetter-Drinfeld module structures
on $V$ and $W$ above.  Observe that the pairing $H^1 \times L^1 \to
k$ gives rise to a pairing $(\, , \, ): V \# H^0 \times W \# L^0 \to
k$, which is explicitly $(v\# x , w\# y) = (v, w)(x , y)$ (by the
definition of $\mathcal R$ and $\mathcal S$). Next, if $\delta: W
\to L^0 \otimes W$ is the coaction, consider the right action
$\leftharpoonup: W \otimes H^0 \to  W$, $w \leftharpoonup a := (a
\otimes Id) \delta(w)$. Then since, by the definition of $W \subset
{\mathcal S}= L^{co \sigma}$, for $w \in W$, we have $\sum w_{(1)}
\otimes w_{(2)}\otimes \sigma (w_{(3)}) =\sum w_{(1)} \otimes
w_{(2)}\otimes 1$, and then
\begin{equation}\label{yd}
(ad (a) v, w) = (v, w\leftharpoonup a), \qquad a\in H^0,
\end{equation}
where $ad(a)v$ denotes the adjoint action, i.e. $ad(a)v = \sum
a_1 vS(a_2).$
 Interchanging r\^oles, we obtain the
relation between the coaction of $V$ and the action of $W$.

{\bf Now  we resume our assumption that  $k$ is algebraically closed
and of characteristic 0, and $\Gamma$ is a finite abelian group.}
Let $V \in {^{\Gamma}_{\Gamma}{\mathcal YD}}$ with basis  $\{ v_1,
\ldots, v_t \}$ where $v_i \in V^{\chi_i}_{g_i}, g_i \in \Gamma,
\chi_i \in \hat{\Gamma}$. Suppose ${\mathcal B}(V)$, the Nichols
algebra, is a finite dimensional Hopf algebra in
$^{\Gamma}_{\Gamma}{\mathcal YD}$. Assume $\chi_i (g_i)$ is a
primitive $r_i$th root of unity for some integer $r_i > 1$, and then
$v_i^{r_i} =0$ in ${\mathcal B}(V)$.


\begin{co}\label{cordual} Let $H = {\mathcal B}(V) \#k [\Gamma]$ for $k,V,\Gamma$ as above.
Then $H^* \cong {\mathcal B} (W) \#k [\hat{\Gamma}]$ where $W=
\oplus_{i=1}^t kw_i \in ^{\hat{\Gamma}}_{\hat{\Gamma}}{\mathcal YD}$
 and $w_i \in W^{g_i }_{\chi_i}$.
\end{co}
\begin{proof} Let $w_i \in W$ be defined by $(v_j,w_i) = \delta_{i,j}$.
Suppose $\delta(w_j) = \sum_{k=1}^t \varphi_{jk} \otimes w_k \in
k[\hat{\Gamma}] \otimes W$. Then for $\leftharpoonup$ the right
action of $\Gamma$ on $W$ described above, $w_j \leftharpoonup h =
\sum_{k=1}^t \varphi_{jk}(h) w_k$ and so by (\ref{yd}),
$\varphi_{jk}(h) = (v_k, w_j \leftharpoonup h) = (ad(h) v_k, w_j) =
\delta_{kj}\chi_k(h)$. Thus $\delta (w_j) = \chi_j \otimes w_j$. A
similar argument shows that for $\chi \in \hat{\Gamma}, \chi \cdot
w_i = \chi(g_i)w_i.$  \end{proof}

Recall from Section 1, that  a pointed Hopf algebra $A$ is called a
lifting of ${\mathcal B}(V) \#k[\Gamma]$ if $gr A$, the associated
graded Hopf algebra, is isomorphic to ${\mathcal B}(V) \#k[\Gamma]$.
The lifting is nontrivial if $A \not \cong gr A$. Now we show that,
although the dual of $H = {\mathcal B}(V)\# k [\Gamma]$ is pointed,
the dual of a nontrivial lifting $A$ of $H$ is not pointed if the
coalgebra projection $\pi$ from $A$ to $A_0$ is not an algebra map.
This result has been shown for various particular cases in the
literature, for example, \cite[p.315]{radminqt}, \cite[Theorem
2.5]{gelaki}, \cite[Cor.5.3]{eg}.

\begin{pr} \label{dualgplikes} Let  $A$ be a nontrivial
lifting of $H= {\mathcal B} (V) \#k[\Gamma]$ such that $\pi: A
\longrightarrow A_0$ is not an algebra map. Then $G(A^*)$, the group
of grouplikes of $A^*$, is isomorphic to a proper subgroup of
$\hat{\Gamma}$ and $A^*$ is not pointed.
\end{pr}
\begin{proof} As above, since $A$ is a lifting of $H$, $A$ is generated as
an algebra by  $(1, g_i)$-primitives $x_i$ with $x_i \in
P'_{(1,g_i)}$ and the elements of $\Gamma$.    The $x_i$ are the
liftings of the $v_i \# 1 \in {\mathcal B} (V) \# k[\Gamma]$. Also
$g_i x_i = \chi_i (g_i) x_i g_i$ and since $\chi_i (g_i) \neq 1$,
all grouplike elements of $A^*$ map the $x_i$ to 0. Thus grouplike
elements of $A^*$ are nonzero only on $A_0= k[\Gamma]$ and so all
grouplikes are of the form $\gamma \cdot \pi$ where $\gamma \in
\hat{\Gamma}$.

However not every $\gamma \in \hat{\Gamma}$ is grouplike in $A^*$.
Since $\pi$ is not an algebra map, there exist $x,y \in A$ such
that $\pi (x) \pi(y) \neq \pi (xy)$ in $k[\Gamma]$. But then
there must exist $\gamma \in \hat{\Gamma}$ such that $\gamma
(\pi(x) \pi(y)) = \gamma \pi (x) \gamma \pi (y) \neq \gamma \pi
(xy),$ for otherwise $\varphi (\pi (x) \pi (y)) = \varphi
(\pi(xy))$ for all $\varphi \in k[\Gamma]^*$. Then $\gamma \pi$
is not an algebra map from $A$ to $k$ and $\gamma \not \in
G(A^*).$ Thus $G(A^*)$ is a proper subgroup of $ \hat{\Gamma}$.

Let $Jac(\cdot)$ denote the Jacobson radical. Since $Jac (k
[\Gamma]) =0$ and $k[\Gamma]$ is commutative, we have that
$k[\Gamma] \cap Jac (A)=0$. For if $z \in k [\Gamma] \cap Jac (A)$,
then $z$ is nilpotent by a theorem of Amitsur \cite[4.20]{lam} and
so lies in $Jac(k [\Gamma]).$ Thus $\dim$ Jac$(A) \leq (\dim A) -
|\Gamma|$ as $k$-spaces, and by \cite[5.1.7]{mont}, $\dim corad
(A^*) = \dim (A) - \dim$ (Jac$(A)) \geq |\Gamma|$. Thus $A^*$ cannot
be pointed. \end{proof}

Thus the coradical of the dual of a nontrivial lifting of a
quantum linear space  contains nontrivial matrix coalgebras.  In
the next sections, we will discuss the number and dimensions of
these matrix coalgebras for some particular examples.

We now use Corollary \ref{cordual} to compute the Drinfel'd double
$D(H)$ for a finite dimensional Hopf algebra $H = {\mathcal B}(V)
\#k [\Gamma]$. Suppose, as above, that  $V = \oplus_{i=1}^t k  v_i$
where $v_i \in V^{\chi_i}_{g_i}$.  Also, as above, define $w_i \in
W^{g_i}_{\chi_i} \subset H^*$ by $(v_j, w_i) = \delta_{i,j}$. Then
$w_i (hv_i) = \chi_i (h)$ and $w_i \in P(H^*)_{\epsilon, \chi_i}$.
Define $w'_i = \chi_i^{-1}
* w_i = ( 1 \# \chi_i^{-1})(w_i \# \epsilon)$ in $H^{*}$ .  Then
$w'_i(z) = \chi_i^{-1}(g_i)  $ for $z$ in $H(1)$ of the form
$hv_i$ and 0 otherwise. In $H^{*}$,  $w'_i$ is $(\chi_i^{-1},
\epsilon_H)$ -primitive. Let $W' = \oplus^t_{i=1} kw'_i   $.

As usual, the left (right) action of $H$ on $H^*$ is denoted $h
\rightharpoonup h^* $ ($h^* \leftharpoonup h$) where     $h
\rightharpoonup h^*(z) = h^*(zh) $ ($h^* \leftharpoonup h(z) =
h^*(hz)$).

\begin{te}
\label{double} Let  $H = {\mathcal B}(V) \# k[\Gamma]$, with
$V=\oplus^t_{i=1} kv_i$ and $v_i \in V^{\chi_i}_{g_i}$. Let $D(H)$
be the Drinfel'd double  of $H$.  Then $gr(D(H)) \simeq {\mathcal B}
(U) \# k[\hat{\Gamma} \times \Gamma]$ where $U \in ^{\hat{\Gamma}
\times \Gamma}_{\hat{\Gamma} \times \Gamma} {\mathcal YD}$ is the
direct sum of $W'$ and $V$ with appropriate actions and coactions,
i.e. $U = \oplus^{2t}_{i=1} ku_i$ where  $u_i \in U^{g_i \times
\chi^{-1}_i}_{\chi^{-1}_i \times 1}$ for $i =1, \ldots, t$ and $u_i
\in U^{g^{-1}_i \times \chi_i}_{\epsilon \times g_i}$ for $i = t+1,
\ldots, 2t$. Also $D(H) \not \cong gr(D(H))$.
\end{te}

\begin{proof} In $H^{*cop}$, the element $w'_i$ is $(\epsilon_H ,
\chi_i^{-1} )$-primitive and $\gamma w'_i \gamma^{-1} = g_i^{**}
(\gamma) w'_i$ for $\gamma \in G(H^{*cop})$ so that $w'_i \in
P(H^{*cop})^{g_i}_{\chi_i^{-1}}$. Thus by an argument similar to
that in Corollary \ref{cordual}, $H^{*cop} \cong {\mathcal B}(W') \#
k[\hat{\Gamma}].$ Since, as coalgebras, $D(H) = H^{*cop} \otimes H$,
by \cite[5.1.10]{mont} or \cite{radminqt}, $D(H)$ is pointed with
group of grouplikes $\hat{\Gamma} \times \Gamma$.

In $D(H)= H^{*cop} \Join H$, the elements $\epsilon \Join v_i$
with $v_i \in V$ are  $(\epsilon \Join 1, \epsilon \Join
g_i)$-primitive and the $w'_i \Join 1$ with $w'_i \in W'$ are
$(\epsilon \Join 1, \chi_i^{-1} \Join 1)$-primitive. For $\gamma
\Join h \in G(D(H))$,  using the formula for multiplication in
$D(H)$ given in \cite[10.3.11 (1)]{mont} from \cite{radminqt} ,
we compute
\begin{eqnarray*}
(\gamma \Join h)(\epsilon \Join v_i)(\gamma^{-1} \Join h^{-1})
&=& (\gamma \Join hv_i)(\gamma^{-1} \Join h^{-1})\\
&=& \sum \gamma [(hv_i)_1 \rightharpoonup \gamma^{-1}
\leftharpoonup \overline{S} ((hv_i)_3)] \Join (hv_i)_2 h^{-1}
\end{eqnarray*}
where the left and right actions of $H$ on $H^*$ are the usual
ones. Now, $\triangle^2 (hv_i) = hv_i \otimes h \otimes h + hg_i
\otimes hv_i \otimes h +hg_i \otimes hg_i \otimes hv_i$ and
$\overline{S} (hv_i) = S^{2r_i -1} (hv_i) = -\chi_i (g_ih)
(g_ih)^{-1} v_i.$ Since $hv_i \rightharpoonup \gamma^{-1}$ and
$\gamma^{-1} \leftharpoonup lv_i$ are $0$ for any $h,l \in
\Gamma$, then

\begin{eqnarray*}
\sum \gamma [(hv_i)_1 \rightharpoonup \gamma^{-1} \leftharpoonup
\overline{S} ((hv_i)_3)] \Join (hv_i)_2 h^{-1} & = & \gamma
[hg_i \rightharpoonup \gamma^{-1} \leftharpoonup h^{-1}] \Join hv_i h^{-1}\\
= \gamma^{-1} (g_i) \epsilon \Join \chi_i (h) v_i &=&(g_i^{**-1}
\otimes \chi_i)(\gamma \otimes h) \epsilon \Join v_i.
\end{eqnarray*}
Then  $\epsilon \Join v_i \in P(D(H))_{\epsilon \Join g_i}^{g_i^{-1}
\Join \chi_i}$ where $g_i^{-1}$ means the usual identification of
$Alg(\hat{\Gamma},k)$ with $\Gamma$. Similarly, we see that $w'_i
\Join 1_H \in P(D(H))^{g_i  \Join \chi_i^{-1}} _{\chi_i^{-1} \Join
1}$. Thus $U=W' \oplus V \in ^{\hat{\Gamma}\times
\Gamma}_{\hat{\Gamma}\times \Gamma}{\mathcal YD}$ with actions and
coactions as in the statement of the theorem.

Now, $gr D(H) \cong {\mathcal R} \#k[\hat{\Gamma} \times \Gamma]$ as
Hopf algebras. Then ${\mathcal B}(W' \oplus V) \subseteq \mathcal
R$. If $1 \leq i,j \leq t$, in the braid matrix for $W' \oplus V$,
$b_{i, t+j} =< g_j^{-1} \Join \chi_j,\chi_i^{-1} \Join 1>= \chi_i
(g_j)$ and $
 b_{t+j,i} = <g_i \Join \chi_i^{-1}, \epsilon _H \Join g_j> = \chi _i (g_j)^{-1}$.
 Thus by \cite[Theorem
2.2]{grana}, $\dim {\mathcal B}(W' \oplus V)= \dim {\mathcal B}(W')
\dim {\mathcal B}(V)= \dim \mathcal R$. Thus $gr D(H) \cong
{\mathcal B} (W' \oplus V)\#k[\hat{\Gamma} \times \Gamma]$.

To show that $D(H) \not \cong gr D(H)$, we compute the products of
the $v_i$ and $w'_j$. By \cite[5.1.10]{mont}, $(\epsilon_H \Join
v_i) (w'_j \Join 1_H) = \sum v_{i_{1}} \rightharpoonup w'_j
\leftharpoonup \overline{S} (v_{i_{3}}) \Join v_{i_{2}} = v_i
\rightharpoonup w'_j \Join 1 + g_i \rightharpoonup w'_j \Join v_i
+ g_i \rightharpoonup w'_j \leftharpoonup \overline{S} (v_i) \Join
g_i$. However $v_i \rightharpoonup w'_j$ and $w'_j \leftharpoonup
\overline{S} (v_i)$ are 0 if $i \neq j$ and then $v_i w'_j =
\chi^{-1} _j (g_i) w'_j v_i$. If $i =j$, then $v_i
\rightharpoonup w_{i}' = \chi_i (g_i)^{-1} \epsilon$ and $g_i
\rightharpoonup w_{i}' \leftharpoonup \overline{S} (v_i) =-
\chi_i (g_i)^{-1} \chi_i$, so that
$$v_i w'_i = \chi_i (g_i)^{-1} w_{i}' v_i + \chi_i (g_i)^{-1}
(\epsilon \Join 1- \chi_i \Join g_i).$$

Thus $\pi$ is not an algebra map and $D(H) \not \cong {\mathcal
B}(U) \#k [\hat{\Gamma} \times \Gamma]$ as algebras. \end{proof}

Note that by adjusting the skew-primitives by a scalar in the
above theorem, we may assume that $ u_{t+j}u_j = \chi^{-1}_j (g_j)
u_j u_{t+j} + (\chi_j^{-1}\Join g_j -\epsilon \Join 1)$.  Recall
from \cite{radminqt} that $D(H)$ is always minimal
quasitriangular and that the dual $D(H)^*$ is also
quasitriangular if and only if $H$ and $H^*$ are quasi-triangular.

For $V$ a quantum linear space, the question of when  Hopf algebras
${\mathcal B}(V) \# k[\Gamma]$ , their liftings and their doubles
are quasi-triangular or ribbon has been considered in
\cite{kauffrad}, \cite{radkauff},  \cite{gelaki}, \cite{pvo},
\cite{adrianacomm}, \cite{adrianajalg}.

\begin{ex}\label{exqls} {\rm  If  $V$ is a quantum linear space, then note that
for $V, W'  \in ^{\hat{\Gamma}\times \Gamma}_{\hat{\Gamma}\times
\Gamma}{\mathcal YD}$ as in Theorem \ref{double}, $U = W' \oplus V
\in ^{\hat{\Gamma} \times \Gamma}_ {\hat{\Gamma} \times \Gamma}
{\mathcal YD}$ is also a quantum linear space.
 Thus if $H = {\mathcal B}(V) \#k[\Gamma]$, $A=D(H)$ is the
 nontrivial lifting of ${\mathcal B} (W' \oplus V) \# k [\hat{\Gamma} \times \Gamma]$
 with lifting matrix ${\mathcal A}  =
  \left [ \begin{array}{ll} 0 & D\\I_t&0 \end{array}
\right] $ where   D is the $t \times t$ diagonal matrix with
$d_{ii} = - \chi_i (g_i)$.
 }\qed
\end{ex}

The quasi-triangularity of $D(H)$ can be used to prove the
quasi-triangularity of some Radford biproducts.

\begin{ex} \label{exponent2} {\rm Let $V,W', U, H,A=D(H)$ be as in
Example \ref{exqls}. Furthermore, let $\chi_i(g_i) =-1$ for
all $1 \leq i \leq t$.   Let $J$ be the Hopf ideal of $A$
generated by the skew-primitives $\chi_i^{-1} \Join g_i -
\epsilon \Join 1
 $ where $1 \leq i \leq t$.
Let $\Gamma ' =( \hat{\Gamma} \times \Gamma)/ \Lambda$ where
$\Lambda$ is the subgroup generated by the $\chi_i^{-1} \Join g_i$.
Then $U$ has a natural $\Gamma '$ grading and also a $\Gamma
'$-action since because $\chi_i(g_i)^2 = 1 = \chi_i(g_j)\chi_j(g_i)$
then $<g_j^{-1} \Join \chi_j, \chi_i^{-1} \Join g_i> = <g_j \Join
\chi_j^{-1}, \chi_i^{-1} \Join g_i>=1$ for all $1 \leq i,j \leq t$.
In fact, it is easy to see that $A/J \cong {\mathcal B} (U) \#
k[\Gamma'] $ where $U$ has this $\Gamma '$-action and grading. Thus
this biproduct is minimal quasi triangular since $A$ is. The
$R$-matrices on $L={\mathcal B} (U) \# k[\Gamma'] $ are discussed in
\cite{adrianajalg}, as well as the existence of ribbon elements for
$L$ and $D(L)$.}\qed
\end{ex}

If $\Gamma \cong C_2$, then  ${\mathcal B} (U) \# k[\Gamma']  =
E(2t)$ and this Hopf algebra has been studied by many authors.

\begin{ex} \label{E(n)}{\rm Let $U$ be the quantum linear space
$U= \oplus^{2t}_{i=1} kv_i \in ^{\hat{\Gamma} \times \Gamma
}_{\hat{\Gamma} \times \Gamma}{\mathcal YD}$ with $\Gamma =<g> \cong
C_2$ and $\hat{\Gamma} = < \chi> \cong C_2$ with $v_i \in U^{(g,
\chi)}_{(1,g)}$ for $ 1 \leq i \leq t$, $v_j \in
U^{(g,\chi)}_{(\chi, 1)}$ for $t +1 \leq j \leq 2t$. Let $A$ be the
nontrivial lifting of ${\mathcal B} (U) \# k [\hat{\Gamma} \times
\Gamma]$ with commutation matrix $\left [\begin{array}{cc} 0 & I_t\\
I_t & 0
\end{array} \right ]$. Then $A$ is minimal quasitriangular since
for $V = \oplus^t_{i=1} kv_i$, we have  $A =D ({\mathcal B} ( V )\#
k [\Gamma]) = D(E(t))$ in the notation of \cite{pvo}. Also it is
shown in \cite{pvo} that $A =D(E(t))$ has quasi-ribbon elements if
and only if $t$ is even, in which case it has 4 quasi-ribbon
elements, two of which are ribbon. And, since it is well-known that
$E(t)$ is quasitriangular  and by \cite{bdg} or by Corollary
\ref{cordual}, $E(t)$ is self-dual, then $A^*=D(E(t))^*$ is also
quasitriangular, but is not pointed by Proposition
\ref{dualgplikes}. If $t =1$, then $A$ is the Drinfel'd double of
the Sweedler Hopf algebra. This 16 dimensional pointed Hopf algebra
is discussed in \cite[Section 2, Case 2]{beattie} and in
\cite[p.562]{cdr} as $H_{(5)}$, and in \cite[p.315]{radminqt}, where
the fact that $A^*$ is not pointed is proved directly. The dual of
$D(E(1))$ will be computed in Section 4.}\qed
\end{ex}

\begin{re}{\rm N. Andruskiewitsch has also pointed out
 that Theorem \ref{double} has an
antecedent in Drinfel'd's paper \cite{drin} where $U_q(g)$ was
introduced as a quotient of the quantum double of $U_q(b)$, and that
this material is strongly tied to the construction in \cite{lusztig}
where Lusztig defines $U_q(g)$ via what are now known as Nichols
algebras. Theorem \ref{double} can also be obtained as a corollary
to a construction of a (possibly infinite) Hopf algebra ${\mathcal
D}(V)$ using the notation and arguments of \cite{lusztig}. Here is
an outline of the proof. For $V,W, \Gamma, \hat{\Gamma}$ as above,
define a bilinear form from $({\mathcal B}(W) \# k[\hat{\Gamma}] )
\times ({\mathcal B}(V)\#k[\Gamma])$ to $k$ by $(\xi \#z | x\#g)=
(\xi, x)(z,g)$. Following \cite[7.2.6]{majid},
 form the graded double ${\mathcal D}(V)$. Identify $V$
(respectively $W$) with subspaces of ${\mathcal D}(V)$ and determine
the commutation relations for $v_iw_j$, $gw_j$, $\gamma v_i$, and
$g\gamma$ where $g \in \Gamma, \gamma \in \hat{\Gamma}$ as was done
in the proof of Theorem \ref{double}. Then ${\mathcal D}(V)$ is
isomorphic to $T(W \otimes k[\hat{\Gamma}] \otimes V \otimes
k[\Gamma]) / {\mathcal J}$ for ${\mathcal J}$ the two-sided ideal
generated by the commutation relations.

As well, the coradical filtration of ${\mathcal D}(V)$ is given by

\hspace*{1in}${\mathcal D}(V)_0 =k [\hat{\Gamma} \times
\hat{\Gamma}]$;

\hspace*{1in}${\mathcal D}(V)_1 = {\mathcal D}(V)_0 +(W
\#k[\hat{\Gamma}]) \Join k[\Gamma] +k [\hat{\Gamma}] \Join (V\#k
[\Gamma]);$

\hspace*{1in}${\mathcal D}(V)_n =({\mathcal D}(V)_1)^n$ for $n >1$.}
\qed
\end{re}

\section{Duals of liftings of quantum linear spaces}

Throughout this section let $V= \oplus^{t}_{i=1} kv_i \in ^{\Gamma}
_{\Gamma} {\mathcal YD}$  be a quantum linear space. We study duals
of Hopf algebras $A$ where $A$ is a nontrivial lifting of ${\mathcal
B}(V) \#k[\Gamma]$ with   matrix ${\mathcal A}$ as in Proposition
\ref{liftqls}. For $v_i \in V^{\chi_i}_{g_i},$ suppose $\chi_i
(g_i)$ is a primitive $r_i$-th root of unity. Let $x_i \in A$ be the
lifting of $v_i \#1$ to $A$, i.e. in $A$, $x_i$ is
$(1,g_i)$-primitive and $hx_i = \chi_i (h) x_i h$ for all $h \in
\Gamma$. Then $A$ has a vector space basis $hz$ with $h \in \Gamma$
and $z \in {\mathcal Z} =\{ x_1^{m_1} \ldots x_t^{m_t}| 0 \leq m_i
\leq r_i -1\}$.

For $\gamma \in \hat{\Gamma}$, define $w( \gamma, x_1^{m_1} \ldots
x_t^{m_t}) \in A^*$  to be the map which takes $hx_1^{m_1} \ldots
x_t^{m_t}$ to $\gamma(h)$ and all other basis elements to 0, so that
${\mathcal W} = \{ w (\gamma, z) | \gamma \in \hat{\Gamma}, z \in
{\mathcal Z} \}$ is a vector space basis for $A^*$. As in Section 2,
we denote $w (\chi_i, x_i)$ by $w_i$ and also $w(\gamma,1)$ is
usually just written $\gamma$.

\begin{lm} \label{multiplication} Let $A$ be a lifting of
${\mathcal B}(V) \# k[\Gamma]$ with  $V$ a quantum linear space.
Let $\gamma, \lambda \in \hat{\Gamma}$, $0 \leq m_i \leq r_i -1$ and
$1 \leq i <j  \leq t$.

(i) Then $$w(\gamma, x_i^{m_i}) *w (\lambda , x_{j}^{m_j} ) =
(\gamma \chi^{-m_i}_i) (g_j^{m_j} ) w(\gamma \lambda,
x_i^{m_i}x_j^{m_j} );$$
$$ w (\lambda , x_j^{m_j})  *w (\gamma, x_i^{m_i}) =
 \lambda(g_i^{m_i}) w(\gamma \lambda, x_i^{m_i}x_j^{m_j}).$$
 In particular, $w_iw_j  = w(\chi_i \chi_j, x_ix_j) =
\chi_i(g_j)w_jw_i$ for $i<j$.

ii) $$\gamma
* w(\lambda, x_1^{m_1} \ldots x_t^{m_t}) = \gamma(g_1^{m_1}
\ldots g_t^{m_t}) w ( \lambda \gamma, x_1^{m_1} \ldots
x_t^{m_t});$$
 $$ w(\lambda , x_1^{m_1} \ldots x_t^{m_t}) * \gamma
=  w(\gamma \lambda, x_1^{m_1} \ldots x_t^{m_t}).$$
 In
particular $\gamma  * w_i  = \gamma(g_i)w_i*\gamma = \gamma (g_i)
w(\chi_i \gamma, x_i)$.

iii) For $0 \leq m_i \leq r_i -1$ we have that
$$ w(\lambda, x_i)*w(\gamma, x_i^{m_i}) =
\lambda(g_i^{m_i})(q^{m_i} + q^{m_i - 1} + \ldots 1)w(\lambda
\gamma, x_i^{m_i + 1})$$ where $q = \chi_i(g_i)^{-1}$, a
primitive $r_i$th root of 1.

The multiplication formulas above show that $\gamma \in
\hat{\Gamma}$ and the elements $w_i$ generate $A^* $ as an algebra
with $\gamma
* w_i = \gamma (g_i) w_i * \gamma$, $w_i * w_j = \chi_i (g_j) w_j
* w_i$ for $i \neq j$, and $w_i^{r_i} =0$.
\end{lm}
\begin{proof}  We prove i). The proofs of ii) and iii) are similar and
straightforward. The map $w(\gamma, x_i^{m_i}) \otimes w (\lambda ,
x_j^{m_j}) $ is nonzero only on elements of $A \otimes A$ of the
form $h x_i^{m_i} \otimes l x_j^{m_j}$. Thus the only basis element
on which $w(\gamma, x_i^{m_i}) *w (\lambda , x_j^{m_j})$ or  $w
(\lambda , x_j^{m_j})
* w(\gamma, x_i^{m_i})$ is nonzero is $h x_i^{m_i}x_j^{m_j}$. Then we have
that
\begin{eqnarray*}
w(\gamma, x_i^{m_i}) *w (\lambda , x_j^{m_j})(h
x_i^{m_i}x_j^{m_j})  &=& w(\gamma,
x_i^{m_i})(hx_i^{m_{i}}g_j^{m_j}) w (\lambda ,
x_j^{m_j})(hx_j^{m_j} )\\
&=& (\gamma \chi_i^{-m_i})(g_j^{m_j} )
 \gamma \lambda (h).
\end{eqnarray*}
Also \begin{eqnarray*}  w (\lambda , x_j^{m_j})
* w(\gamma, x_i^{m_i})(h x_i^{m_i}x_j^{m_j})
&=& w (\lambda , x_j^{m_j})(h g_i^{m_i} x_j^{m_j})w(\gamma,
x_i^{m_i})(hx_i^{m_i}) \\
&=& (\lambda \gamma)(h) \lambda(g_i^{m_i}),
\end{eqnarray*}
which proves i).
 \end{proof}

\begin{pr}\label{basic}
Let $V= \oplus^t_{i=1} kv_i$ be a quantum linear space,
  and let
$A$ be a  lifting of ${\mathcal B}(V) \#k[\Gamma]$  with matrix
${\mathcal A} $ such that for $i \neq j$, if $\alpha_{i,j} \neq 0$,
then $\chi_i = \chi_j$ so that $\chi_i^2 = \chi_j^2 = \epsilon$ and
$r_i = r_j =2$. We make the following definitions:
\begin{center}
\begin{description}
 \item[ $\Theta$] is the subgroup of $\hat{\Gamma}$ generated by the
$\chi_i, 1 \leq i \leq t$;
\item [ $\Omega$] is the subgroup of $\hat{\Gamma}$ defined by
$\Omega = \{ \gamma \in \hat{\Gamma}: \gamma
(\alpha_{i,i}(g_i^{r_i} - 1)) = 0 = \gamma(\alpha_{i,j}(g_ig_j -
1))  \mbox{ for } 1 \leq i , j \leq t , i \neq j \};$
\item [ $\Gamma'$]  is the subgroup
of $\Gamma$ generated by $\{ g_i^{r_i} |\alpha_{i,i} \neq 0 \}
\cup \{ g_ig_j | \alpha_{i,j} \neq 0 \}.$
\end{description} \end{center}
Let  $\overline{\Gamma} = \Gamma /
\Gamma' $ and note that $\hat{\overline{\Gamma}} \cong \Omega.$
For $g \in \Gamma$, $\overline{g}$ will denote its image in
$\overline{\Gamma}$. For $\gamma \in \hat{\Gamma}$,
$\overline{\gamma} $ will denote its image in
$\hat{\Gamma}/\Theta$.
     Then
 \begin{enumerate}
\item   $\Theta$ is a subgroup of $\Omega$ and $\Omega$ is the  group
of grouplikes in $A^*$.
\item There exist subHopf algebras of $A^*$ isomorphic to
${\mathcal B}(W) \#k [\Theta]$ and ${\mathcal B}(W) \# k [\Omega]$
where $W = \oplus ^t_{i=1} kw_i$ with $w_i \in
W^{\overline{g_i}}_{\chi_i}$. Denote by $C(\overline{1})$ the
pointed subHopf algebra isomorphic to ${\mathcal B}(W) \#k [\Theta]$
and by $B$ the  subHopf algebra isomorphic to ${\mathcal B} ( W) \#
k [\Omega]$.
\item For $\gamma \in \hat{\Gamma}$, let $  C(\overline{\gamma})=
\gamma  C(\overline{1}) =  C(\overline{1})\gamma$ for any $\gamma
\in \overline{\gamma}.$  $C(\overline{\gamma})$ is a subcoalgebra
of $A^*$ and as a coalgebra, $A^* \cong \oplus _{
\overline{\gamma} \in \hat{\Gamma}/\Theta} C(\overline{\gamma}).$
If $\gamma \in \Omega$, then $C(\overline{\gamma}) \cong
C(\overline{1})$ as a coalgebra and $B \cong \oplus _{
\overline{\gamma} \in \Omega/\Theta} C(\overline{\gamma})$ as a
coalgebra. Thus corad$A^* \cong k [\Omega] \oplus
   \oplus_{ \scriptstyle \overline{\gamma} \in
\hat{\Gamma}/ \Theta \atop \scriptstyle \gamma \not \in \Omega}
(corad C(\overline{\gamma})).$
  \item For all
  $\overline{\gamma}, \overline{\lambda} \in
  \hat{\Gamma}/\Theta$, $C(\overline{\gamma}) C(\overline{\lambda})
  = C(\overline{\gamma \lambda})$ and also
  $S(C(\overline{\gamma})) = C(\overline{\gamma^{-1}})$.
\end{enumerate}
\end{pr}

\begin{proof} \par (1)  Clearly $\gamma \in \hat{\Gamma}$ is
grouplike in $A^*$ if and only if $\gamma \in \Omega$. Now, if
$\alpha_{i,i} \neq 0$ then $\chi_i^{r_i} = \epsilon$ and thus for $j
\neq i$, $\chi_j(g_i^{r_i} )= \chi_i^{r_i}(g_j^{-1}) = 1$ so
$\chi_j(\alpha_{i,i}(g_i^{r_i}-1)) = 0$ for all $i$. Also, if
$\alpha_{i,j} \neq 0$, then for $k \neq i,j$, $\chi_k(g_ig_j) =
(\chi_i\chi_j)(g_k) = 1$. As well, $\chi_i(g_ig_j) =
\chi_i(g_i)\chi_i(g_j) = \chi_i(g_i ) \chi_j^{-1}(g_i) =
(\chi_i\chi_j)(g_i) = 1$ since we are assuming here that $\chi_j^2 =
\epsilon$. \par (2) Let $\overline{V} = \oplus _{i=1}^t k
\overline{v_i} \in ^{\overline{\Gamma}}_{\overline{\Gamma}}{\mathcal
YD}$ with $\delta (\overline{v_i}) = \overline{g_i} \otimes
\overline{v_i}, \overline{h} \to \overline{v_i} = \overline{h \to
v_i} = \chi_i (h) \overline{v_i}$ where $\overline{g_i},
\overline{h}$, denote the images of $g,h$ in $\overline{\Gamma}$.
The action of $\overline{\Gamma}$ on $\overline{V}$ is well-defined
by the arguments in (1). There is an exact sequence of Hopf
algebras,

$$\xymatrix{1 \ar[r]& k[\Gamma']
\ar[r]&A \ar[r]^-{\phi} &{\mathcal B}(\overline{V}) \# k
 [\overline{\Gamma}]  \ar[r] & 1}$$ where $\phi$ is defined by $\phi
(hx_1^{k_1}x_2^{k_2}\ldots x_t^{k_t}) = \overline{h}
\overline{v_1}^{k_1}\overline{v_2}^{k_2}\ldots \overline{v_t}^{k_t}
$. By Corollary \ref{cordual} the dual of ${\mathcal B}
(\overline{V}) \#k \overline{[\Gamma]}$ is ${\mathcal B} ( W) \# k
[\Omega]$ where $W = \oplus^t_{i=1} kw_i$ with $w_i \in
W^{\overline{g_i}}_{\chi_i}$ and so ${\mathcal B}(W) \#k [\Omega]$
is isomorphic to a pointed sub Hopf algebra of $A^*$. We denote this
subHopf algebra by $B$ and note that $C(\overline{1}) = {\mathcal
B}(W) \#k [\Theta]$ is a subHopf algebra of $B$.

\par(3) Since $k[\Theta] \subseteq C(\overline{1})$, a subHopf algebra of
$A^*$, $C(\overline{\gamma})$ is well-defined and  Lemma
\ref{multiplication}  implies  $ \gamma C(\overline{1})=
C(\overline{1})\gamma.$ If $\gamma$ is grouplike, then the map $z
\to \gamma z$ from $C(\overline{1})$ to $\gamma C(\overline{1})$ is
a coalgebra isomorphism. If $\gamma \not \in \Omega$, then
$\triangle (\gamma) = \gamma \otimes \gamma + z(\gamma) \in A^*
\otimes A^*$. Here $z(\gamma)$ is a sum of terms of the form $\xi
w(\lambda, x_1^{n_1} \ldots x_t^{n_t} ) \otimes w(\lambda
',x_1^{m_1} \ldots x_t^{m_t})$ where we have $\xi \lambda(h) \lambda
'(l) = \gamma(h  x_1^{n_1} \ldots x_t^{n_t} l  x_1^{m_1} \ldots
x_t^{m_t}) =    \lambda ''(l)\gamma ( hl x_1^{n_1} \ldots
x_t^{n_t}x_1^{m_1} \ldots x_t^{m_t})$ where $\lambda '' \in \Theta$.
Since this holds for all $h,l$ we must have $\xi = \gamma( x_1^{n_1}
\ldots x_t^{n_t} x_1^{m_1} \ldots x_t^{m_t})$ and $\lambda =
\gamma$, $\lambda ' = \gamma \lambda ''$ with $\lambda '' \in
\Theta$. Thus $w(\lambda, x_1^{n_1} \ldots x_t^{n_t}) \otimes
w(\lambda ',  x_1^{m_1} \ldots x_t^{m_t}) \in C(\overline{1})\gamma
\otimes C(\overline{1})\gamma$ and $C(\overline{\gamma}) $ is a
coalgebra as required. Clearly $C(\overline{1}) \gamma $  has $k$
basis ${\mathcal W}_\gamma = \{ w(\gamma \chi, z)| \chi \in \Theta,
z \in {\mathcal Z}\}$ so that $A^* = \oplus_{\overline{\gamma} \in
\hat{\Gamma}/\Theta} C(\overline{\gamma}) $ as required.

\par(4)
Since $C(\overline{1})$ is a subHopf algebra,
$$C(\overline{\gamma})C(\overline{\lambda}) = \gamma  C
(\overline{1}) C(\overline{1}) \lambda = \gamma C(\overline{1})
\lambda = (\gamma \lambda)C(\overline{1}) =C(\overline{\gamma
 \lambda}).$$
 If $\gamma \in \hat{\Gamma}$,  since $S_{A^*}(\gamma) =
 (S_A)^*(\gamma) = \gamma^{-1}$, the last statement is clear.

  \end{proof}

We describe $A^*$ when $\dim V=1$. It is shown in \cite[2.3]{rad}
using the algebra structure of $A$ that for $\Gamma = <g>$ cyclic of
order $np$ and $V=kv=V_g^{\chi^p}$where $\chi$ generates
$\hat{\Gamma}$, then for $A$ the unique nontrivial lifting of
${\mathcal B}(V) \# k[\Gamma]$, the dual  $A^*$ contains $p -1$
matrix coalgebras of dimension $n^2$. Also it is shown in
\cite{gelaki} that $A^*$ contains a nontrivial matrix subcoalgebra.
 We generalize this description to
$\Gamma$ any finite abelian group and explicitly calculate   the
matrix coalgebras.

\begin{te} \label{rank1} Let $V = kv \in {^{\Gamma}_{\Gamma}{\mathcal YD}}$
with $v \in V_g^\chi$ and   $\chi(g)$   a primitive $r$th root of
unity. Suppose $\chi^r =\epsilon, g^r \neq 1$ and $A$ is the
(only)
 nontrivial lifting
of ${\mathcal B}(V) \# k[\Gamma]$.   Then as a coalgebra, $A^* \cong
\oplus_{\overline{\gamma} \in \hat{\Gamma}/\Theta}
C(\overline{\gamma})$ where $C(\overline{\gamma})\cong {\mathcal
M}^c(r,k)$  if $\gamma \not \in \Omega$ and
$C(\overline{\gamma})\cong C( \overline{1}) \cong {\mathcal B}(W) \#
k [\Theta]$ as coalgebras with $W = k \cdot w, w \in
W_\chi^{\overline{g}}$ if $\gamma \in \Omega$. Thus corad $(A^*)$ is
the sum of the group algebra $k [\Omega]$ and $m $ $r \times r$
matrix coalgebras where
  $\displaystyle mr =
|\hat{\Gamma}| - |\Omega|.$ As an algebra $A^*$ is generated by
$\hat{\Gamma}$ and $w$ where $\gamma *w = \gamma (g) w* \gamma$
for $\gamma \in \hat{\Gamma}$  and $w^r =0.$
\end{te}
\begin{proof} Let $\gamma \in \hat{\Gamma}- \Omega$ and we show that
$C(\overline{\gamma})$ is an $r \times r$ matrix coalgebra by
explicitly exhibiting a matrix coalgebra basis. For $1 \leq i,j \leq
r$, define $e(i,j) = w(\gamma \chi^{i-1}, x^{i-j})$ if $i \geq j$
and $e(i,j) = \eta w(\gamma \chi^{i-1}, x^{r+i-j})$ if $i <j$ where
$\eta = \gamma (g^r-1) \neq 0$. Then we claim that $\triangle
(e(j,k)) =( m^*_A)(e(j,k))$ is $\sum^r_{s=1} e(j,s) \otimes e(s,k).$
First we note that for $j \geq k$,
\begin{eqnarray*}
\triangle w(\gamma \chi^{j-1}, x^{j-k})&=& \sum^{k-1}_{s=1}
w(\gamma \chi^{j-1}, x^{j-s}) \otimes \eta w(\gamma \chi^{s-1},
x^{r+s-k}) \\
&+& \sum^j_{s=k} w(\gamma \chi^{j-1}, x^{j-s}) \otimes w (\gamma
\chi^{s-1}, x^{s-k})\\
&+& \sum^r_{s=j+1} \eta w (\gamma \chi^{j-1}, x^{r+j-s}) \otimes w
(\gamma \chi^{s-1}, x^{s-k}).
\end{eqnarray*}
Now in the first sum, $s< k \leq j$ and so we have $e(j,s) =
w(\gamma \chi^{j-1}, x^{j-s})$ and $e(s,k) = \eta w(\gamma
\chi^{s-1}, x^{r+s-k}),$ and the sum is $\sum^{k-1}_{s=1} e(j,s)
\otimes e(s,k)$. The others are similar and we conclude that
$\triangle (w(\gamma \chi^{j-1}, x^{j-k})) = \triangle (e(j,k)) =
\sum^r_{s=1} e(j,s) \otimes e(s,k).$

Similarly, if $j<k$,
\begin{eqnarray*}
\triangle \eta w(\gamma \chi^{j-1}, x^{r+ j-k})&=& \eta [
\sum^{j}_{s=1} w(\gamma \chi^{j-1}, x^{j-s}) \otimes w(\gamma
\chi^{s-1},
x^{r+s-k})] \\
&+& \eta [\eta \sum^{k - 1}_{s=j+1} w(\gamma \chi^{j-1}, x^{j-s})
\otimes w (\gamma
\chi^{s-1}, x^{r + s-k})] \\
&+& \eta [\sum^r_{s=k}  (\gamma \chi^{j-1}, x^{r+j-s}) \otimes w
(\gamma \chi^{s-1}, x^{s-k})].
\end{eqnarray*}

Note that the middle sum is empty if $j + 1 =k$. It is easily
verified that this sum is $\sum_{s=1}^r e(j,s) \otimes e(s,k)$.
  The last statement follows from Lemma
\ref{multiplication}. \end{proof}

\begin{re}
\label{Jac} For $A,V, \Gamma, g, \chi$ as in Theorem \ref{rank1},
Jac$(A)$ is the ideal $I$ generated by $(1+ g^r + \ldots
+g^{(s-1)r}) x$ where $s$ is the order of $g^r$. For $I$ is
clearly contained in Jac$(A)$ since $I$ is a nil ideal and $\dim
I = | \Gamma/<g^r>| (r-1) = |\Omega|(r-1)$, while $\dim corad
(A^*) = |\Omega| + r(|\Gamma| -|\Omega|) = \dim A- |\Omega|(r-1).$
\qed \end{re}

\begin{co}
\label{quotient} Let $A$ be a lifting of ${\mathcal B}(V)
\#k[\Gamma]$ where $V = \oplus^t_{i=1} kv_i$ is a quantum linear
space and $x_i$ is the lifting of $v_i \#1$ to $A$. Suppose that the
matrix ${\mathcal A}$ has $\alpha_{1,1} = 1$  and has all other
entries $0$, i.e. $x_i x_j = \chi_j (g_i) x_jx_i,$ for $1 \leq i,j
\leq t$, and $  x_j^{r_j} =0$ for $j= 2, \ldots, t$ while $x_1^{r_1}
= g_1^{r_1} -1$.
  Then $A^*$ has coradical $C \cong k
[\Omega] \oplus \oplus_{ \scriptstyle \overline{\gamma} \in
\hat{\Gamma}/ \Theta \atop \scriptstyle \gamma \not \in \Omega}
{\mathcal M}^c (r_1,k)$.
\end{co}

\begin{proof} Let $I$ be the Hopf ideal of $A$ generated by the
skew-primitives $x_2, \ldots, x_t$. Then there is a Hopf algebra map
from $A$ onto $A_1 = A/I$, and $A_1$ is isomorphic to the lifting of
${\mathcal B}(kv_1) \# k [\Gamma]$ considered in Theorem
\ref{rank1}. Thus $A^*$ has a sub Hopf algebra isomorphic to $A_1^*$
and by Theorem \ref{rank1}, the coradical of $A_1^*$ is isomorphic
to $C$. Now, Jac$(A)$ contains the nil ideals $N_1 = <x_2, \ldots,
x_t>$ and $N_2 =<(1+ g_1^{r_1} + \ldots + g_1^{(m_1-1)r_1}) x_1>$
where $m_1r_1$ is the order of $g_1$. Thus $\dim Jac(A) \geq
|\Gamma| r_1 (r_2 \ldots r_t-1) +|\Omega| (r_1-1) r_2 \ldots r_t -
|\Omega| (r_1 -1) (r_2 \ldots r_t -1)= |\Gamma| r_1 (r_2 \ldots
r_t-1) +|\Omega| (r_1-1),$ i.e. Jac$(A)$ has dimension greater than
or equal to $\dim N_1+\dim N_2-\dim(N_1 \cap N_2)$. Thus $\dim
(corad A^*) \leq |\Gamma| r_1r_2 \ldots r_t - |\Gamma| r_1 (r_2
\ldots r_t-1) -|\Omega| (r_1-1) $ $= |\Gamma| r_1- |\Omega| (r_1-1)
= |\Omega| +( |\Gamma| - |\Omega|)r_1 = \dim C$. \end{proof}

Now we consider quantum linear spaces of dimension 2.  The next
result is proved in \cite{rad}   for the group $\Gamma$ cyclic
but without our restrictions on $r$.  The approach is rather
different.

\begin{te}\label{rank2} Let $V= kv_1 \oplus kv_2$ be a quantum
linear space with $v_1  \in V^\chi_g$ and $v_2 \in
V_{g^n}^{\chi^m}$ and such that

i) $\chi(g)$ and $\chi^m(g^n)$ are primitive $r$th roots of unity
for $r$ an  odd squarefree integer.

ii) $\chi$ and $\chi^m$ have order $r$ in the group
$\hat{\Gamma}$, i.e. $(m,r)=1 =(n,r).$

iii) $g^r \neq 1$ and $g^{nr} \neq 1$.

Let $A$ be the lifting of ${\mathcal B}(V) \# k[\Gamma]$ with matrix
${\mathcal A} =I_2$ . Then for $\gamma \not \in   \Omega  = G(A^*)$,
$C(\overline{\gamma}) \cong \oplus^r_{i=1} {\mathcal M}^c (r,k)$,
and $corad (A) \cong k [\Omega] \oplus \oplus_{ \scriptstyle
\overline{\gamma} \in \hat{\Gamma}/ \Theta \atop \scriptstyle \gamma
\not \in \Omega} ( \oplus^r_{i=1} {\mathcal M}^c (r,k)).$
\end{te}
\begin{proof} Let $x$ and $y$ be the liftings of $v_1$ and $v_2$ to $A$.
Then $yx =q^n xy$ where $q = \chi (g)$ is a primitive $r$th root of
unity and $x^r = g^r -1$, $y^r = g^{nr} -1$. The dimension of
$C(\overline{\gamma})$ for any $\overline{\gamma}$ is $r^3$.

Let $\gamma \not \in \Omega = G(A^*).$ Since $C(\overline{\gamma})$
has no grouplike elements, $C(\overline{\gamma})$ must contain a
matrix coalgebra $M \cong {\mathcal M}^c (s,k)$ for some integer $s
\geq 2.$ Thus $M$ has a matrix coalgebra basis $e(i,j), 1 \leq i,j
\leq s$ with $\triangle (e(i,i)) = e(i,i) \otimes e(i,i) + z_i \in
C(\overline{\gamma}) \otimes C(\overline{\gamma})$ and $\epsilon
(e(i,i)) =1$.

The only elements $z $ in the basis ${\mathcal W}_\gamma =\{
w(\gamma \chi^i, x^j y^k)$, $0 \leq i,j,k \leq r-1\}$ for
$C(\overline{\gamma})$, such that a scalar multiple of $z \otimes z$
is a summand of $\triangle w$ for some basis element $w \in
{\mathcal W}_{\gamma}$ are those of the form $ z = w(\gamma \chi^i,
x^jy^k)$ where $\chi^j\chi^{mk} = \epsilon$, i.e., $x^jy^k$ commutes
with elements of $\Gamma$. Thus for every $l$, $1 \leq l \leq s$,
then $ e(l,l)$ is a linear combination of such elements. If $e(l,l)
= \sum^{r-1}_{i=0}a_i \gamma \chi^i + w$, for $w$ a linear
combination of basis elements from ${\mathcal W}_\gamma \setminus
\hat{\Gamma}$, then it is clear that exactly one $a_i$ is 1 and the
rest are 0 since terms of the form $\gamma \chi^i \otimes \gamma
\chi^j$ with $i \neq j$ do not occur in $\triangle (z)$ for any $z
\in A^*$ . For $w =w (\gamma \chi^i, x^jy^k) \in {\mathcal
W}_\gamma$, $\triangle w= \gamma \chi^i \otimes w + w \otimes \gamma
\chi^i \chi^{-j} \chi^{-km} + w'$ where $w'$ is a linear combination
of terms from  $ {\mathcal W}_\gamma \setminus \hat{\Gamma} \otimes
{\mathcal W}_\gamma \setminus \hat{\Gamma}$. Thus $e(l,l) = \lambda
+ \sum^{r-1}_{i=1} a_i w(\lambda, x^{nr-im} y^i)=\lambda +
\sum^{r-1}_{i=1} a_i z_{i,\lambda}$, for some $\lambda \in
\overline{\gamma}$, where $z_{i,\lambda} = w(\lambda,x^{nr-im} y^i)
$ with $n$  the integer such that $0 \leq nr -im \leq r-1$. (This
also shows that the $e(i,j), i \neq j, $ are linear combinations of
elements of ${\mathcal W}_\gamma \setminus \hat{\Gamma}.$) Now,
 in the expression for $\triangle (\lambda)$, a nonzero scalar multiple of
 $w(\lambda, x^j) \otimes  w(\lambda \chi^{-j}, x^{r-j}) $
 occurs for all $1 \leq j \leq r-1$.  Also no such summand occurs in
 $\triangle z_{i,\lambda}$ for any $i$ and so $w(\lambda, x^j)$ must
 be a summand in some $e(l,m)$.  But then
$e(m,m)$ must contain $\lambda \chi^{-j}.$  Then, for all $\theta
\in \overline{\gamma}$, $\theta $ is a summand of some matrix
coalgebra basis element $e(m,m)$. Thus the dimension of $M$ is at
least $r^2$.

Since $\chi^i$ is grouplike,  the map $z \to \chi^i z \chi^{-i}$
from $C( \overline{\gamma})$ to $C( \overline{\gamma})$ is a
coalgebra isomorphism. Suppose the $r$ subcoalgebras $\chi^i M
\chi^{-i}$ of $C( \overline{\gamma})$, $i=0, 1, \ldots, r-1$, are
not all distinct (this will be the case if $s>r$  ) and so
$\chi^i M \chi^{-i} =M$ for some $i$. Let $\Theta' = \{\chi^i |
\chi^iM\chi^{-i} = M \} \subseteq <\chi>$, and let $k$ to be the
smallest positive integer such that $\chi^k \in \Theta'$. Then
$\Theta ' = <\chi^k>$ and $r=kt$ for some $t$.

For $\theta, \lambda \in C(\overline{\gamma})$, denote by
$E(\theta, \lambda)$ the subspace of $M$ spanned by those
$e(i,j)$ such that $\triangle e(i,j)$ contains $\theta \otimes
e(i,j)$ and $e(i,j) \otimes \lambda$ as summands. By the previous
discussion, $\dim E (\theta, \lambda) \geq 1$ for all $\theta,
\lambda \in \overline{\gamma}$. Note that $\dim E (\theta,
\theta)$ is a square, since if $e(i,i), e(j,j) \in E(\theta,
\theta)$ then $e(i,j), e(j,i) \in E(\theta, \theta).$ Conversely,
if $e(i,j) \in E(\theta,\theta)$, then $e(i,i)$ and $e(j,j)$ lie
in $E(\theta,\theta).$ Also if $\dim E(\theta, \theta) = l^2$ and
$\dim E(\lambda, \lambda) =m^2$ then $\dim E(\theta, \lambda) =
\dim E(\lambda, \theta) =lm.$ This is because if $e(i,i) \in
E(\theta, \theta)$ and $e(j,j) \in E(\lambda, \lambda) $ then
$e(i,j) \in E(\theta, \lambda)$ and $e(j,i) \in E(\lambda,
\theta)$ and conversely. Note that $\chi^k E(\theta, \lambda)
\chi^{-k} = E(\theta, \lambda)$ for all $\theta, \lambda \in
\overline{\gamma}.$

Let $E(\theta, \theta)$ contain a vector $v$ such that $v= a_0
\theta + \sum^{r-1}_{i=1} a_iz_{i,\theta}$ for $z_{i,\theta}
=w(\theta, x^{nr-im}y^i)$ as above, with $a_i \neq 0$ for $0 \leq
i \leq r-1$. Now $\chi^k v \chi^{-k} =a_0 \theta +
\sum^{r-1}_{i=1}
 \chi^k (g^{nr-im} g^{in})
a_i z_i$ where $\chi^k (g^{nr-im} g^{in}) = q^{ki(m-n)}$. Since
$V$ is a quantum linear space, $q^{m+n}=1$ and so $m+n =pr$ for
some $p$ and $m-n=2m-pr$. Thus $q^{ki(m-n)} = (q^{2m})^{ki}$ and
since $(2,r) =(m,r) =1,$ $q^{2m}$ is a primitive $r$th root of
unity, and  $q'=q^{2mk}$ is a primitive $t$th root of unity. Thus
the $t$ vectors $\chi^{ik} v \chi^{-ik}$, $0 \leq i \leq t-1$, are
linearly independent since the rank of the Vandermonde matrix
$V(1, q', q'^{2}, \ldots , q'^{t-1})$ is $t$. So $\dim E(\theta,
\theta) \geq t$, and since if $t
>1, t$ is not a square by the assumption that $r$ is square-free, $\dim E(\theta, \theta)
>t$. Thus, if all $E (\theta, \theta)$ contain such a vector $v$,
then $\dim M = \sum_{\theta, \lambda \in \overline{\gamma}} \dim
E(\theta, \lambda) > r^2 t$. Since the isomorphic coalgebras $M,
\chi M \chi^{-1}, \ldots, \chi^{k-1} M \chi ^{- (k-1)}$ are distinct
by our choice of $k$, $\dim C(\overline{\gamma}) > r^2 tk =r^3$,
which is a contradiction. Thus $k=r,t=1, $ $M$ has dimension $r^2$
and $C(\overline{\gamma}) \cong \oplus^r_{j=1} {\mathcal M}^c
(r,k)$.

It remains to show that every $E(\theta,\theta)$ contains a
vector $v$ of the required form.

Let $e(l,l) \in E(\theta, \theta)$ with $\dim E(\theta, \theta)
>1$ and suppose $a_1=0$ in the expression for $e(l,l)$. Since
$\triangle \theta$ is a summand of $\triangle e(l,l), \zeta
z_{1,\theta} \otimes z_{r-1, \theta}$ is a summand of $\triangle
e(l,l)$, and $\zeta' z_{1,\theta}$ is a summand of $e(l,j)$ for
some $j \neq l$, with $e(l,j), e(j,j), e(j,l) \in E (\theta,
\theta)$. (Here $\zeta, \zeta'$, etc. denote nonzero scalars.)
Similarly, a nonzero multiple of $z_{2,\theta}$ is a summand in
$e(l,l)$ or in some $e(l,j)\in E(\theta,\theta)$.  Continuing in
this way we find a vector $v$ in $E(\theta, \theta)$ as required,
by taking some linear combination of $e(l,l)$ and the $e(l,j)$.
\end{proof}

\begin{re}\label{requal2}
In the setting of Theorem \ref{rank2}, if $r=2$, then $q^{2nki} =
1$ and   $\chi^ke(\lambda, \lambda) \chi^{-k} = e(\lambda,
\lambda)$ so that $M$ is invariant under the action of $\chi.$
\qed
\end{re}

\begin{co} \label{Jacobsonrank2} For $A$ as in Theorem \ref{rank2},
$I = \{ h (1 + g^r + \ldots + g^{r(l-1)}) x^i y^j | h \in
\Gamma$, $l$ is the order of $g^r$ in $\Gamma$, $0 \leq i,j \leq
r-1$ and $i + j >0 \}$ is the Jacobson radical of $A$.
\end{co}

\begin{proof} Clearly $I $ is a nil ideal in $A$ and has dimension
$|\Omega|(r^2 -1)$. But
 $$
\dim Jac(A) = | \Gamma | r^2 - \dim corad (A^*) = | \Gamma | r^2
-( |\Omega| + \left ( \frac{|\Gamma | - |\Omega|}{r}   \right )
r^3) = |\Omega|(r^2 -1).$$
  \end{proof}

\begin{co} \label{quotient2}
Let $A$ be a lifting of ${\mathcal B}(V) \#k[\Gamma]$ where $V=
\oplus^t_{i=1}kv_i$ is a quantum linear space such that $V'= kv_1
\oplus kv_2$ satisfies the conditions in Theorem \ref{rank2}, and
$A$ has lifting matrix ${\mathcal A} = \left [ \begin{array}{c|c}
I_2 &0
\\ \hline 0&0 \end{array} \right ]$. Then in $A^* , corad
C(\overline{\gamma}) \cong \oplus^r_{i=1} {\mathcal M}^c (r,k)$ for
$\gamma \not \in \Omega$.
\end{co}

\begin{proof} Let $N_1 =\{ h(1+ g^r + \ldots + g^{r(l-1)} ) x_1^{m_1}
x_2^{m_2} \ldots x_t^{m_t} | h \in \Gamma$, $l$ is the order of
$g^r$, $0 \leq m_i \leq r_i -1$ and $m_1 + m_2 >0 \}$ and let $N_2 =
\{ hx_1^{m_1} x_2^{m_2} x_3^{m_3} \ldots x_t^{m_t} | h \in \Gamma, 0
\leq m_i \leq r_i -1$, $m_3 + \ldots + m_t >0 \}$. Since $N_1$ and
$N_2$ are nil ideals of $A$, then $Jac(A)$ has dimension at least
$|\Omega|(r^2 -1) r_3 \ldots r_t + | \Gamma | r^2 (r_3 \ldots r_t
-1) - |\Omega|(r^2 -1) (r_3 \ldots r_t-1)= |\Gamma| r^2 (r_3 \ldots
r_t-1) + |\Omega|(r^2-1).$ Thus $\dim corad (A^*) \leq | \Gamma |
r^2 r_3 \ldots r_t - | \Gamma | r^2 (r_3 \ldots r_t-1) -|\Omega|
(r^2 -1) = (|\Gamma | -|\Omega|) r^2 +|\Omega|$. Now, as in the
proof of Corollary \ref{quotient}, let $I$ be the Hopf ideal of $A$
generated by the skew-primitives $x_3, \ldots, x_t$. Then $A_1 =
A/I$ is isomorphic to the lifting of ${\mathcal B}(kv_1 \oplus kv_2)
\#k[\Gamma]$ with matrix ${\mathcal A} =I_2$ and so $A_1^* \subseteq
A^*$ has coradical isomorphic to $k [\Omega] \oplus \oplus_{
\scriptstyle \overline{\gamma} \in \hat{\Gamma}/ \Theta \atop
\scriptstyle \gamma \not \in \Omega} (\oplus^r_{i=1} {\mathcal M}^c
(r,k))$ by Theorem \ref{rank2}. But the dimension of this space is
$|\Omega| + ( \frac{|\Gamma|-|\Omega|}{r} ) r^3 = |\Omega| + ( |
\Gamma| -|\Omega|) r^2$. \end{proof}

In the next section we do some detailed computations of bases for
matrix coalgebras in the duals of some liftings of quantum linear
spaces.

\section{Examples}
\subsection{Duals with grouplikes of order 3.}
 The   two examples in this subsection illustrate Theorem
 \ref{rank2}, i.e. they are duals of liftings of quantum linear
 spaces with matrix ${\mathcal A} =I_2$.

\begin{re}\label{choosee}{\rm
In the setting of Theorem \ref{rank2}, in order to choose
$e(i,i)$ in a matrix coalgebra basis, each summand $\alpha w$ of
$e(i,i)$ must be one such that $ \alpha w \otimes \alpha w$
 can occur as a summand in $\triangle(\beta w')$ where $w, w'$ lie in a
 basis $\mathcal W$ for $A^*$.
 This was noted in the proof of Theorem \ref{rank2}.

\par Suppose $v_1 \in V_g^{\chi}$ and $v_2 \in V_g ^ {\chi^{-1}}$ so
that in $C(\overline{\gamma})$, $e(1,1) =
\sum_{i=0}^{r-1}a_iw(\gamma, (xy)^i).$ Then $a_0 = 1.$ Clearly $a_2
= a_1^2, a_4 = a_2^2$, etc. If $r=3$, then $a_2 = a_1^2$ and $
\eta^2 a_1 = a_2^2$ where $\eta = \gamma(g^3-1)$. Then $\eta^2a_1 =
a_1^4$ so $a_1^3 = \eta^2$ and $a_1 = q \eta^{2/3}$ where $q^3 =1.$
Similar computations can be done for other values of $r.$ For
instance, if $r=5$, then $a_2= a_1^2, a_4 = a_2^2, \eta^2 a_3 =
a^2_4, \eta^2 a_1 = a_3^2$ where here $\eta = \gamma(g^5-1).$  Then
$\eta^6 a_1 = \eta^4a_3^2 = a_4^4 = a_2^8 = a_1^{16}$ and thus
$a_1^{15} = \eta^6$ so if $a_1 = q \eta ^{6/15} = a \eta ^{2/5}$
where $q^{5}=1$ then $a_i = q ^i\eta^{2i/5}.$ \qed }
\end{re}

\begin{ex}\label{order3I} {\rm For $\Gamma$ finite abelian as
usual, let $V = kv_1 \oplus kv_2$ with $v_1 \in V_g^\chi, v_2 \in
V_g^{\chi^{-1}}$ where $\chi^3 = \epsilon$ and $g^3 \neq 1$. Let $A$
be the nontrivial lifting of ${\mathcal B}(V) \# k[\Gamma]$ with
matrix ${\mathcal A} =  I_2$. Denote by $x,y$ the liftings of the
$v_i \#1$ to $A$ and write $q = \chi(g)$, a primitive cube root of
unity. Then $yx = qxy$, $x^3 =y^3 = g^3 -1$.  By Proposition
\ref{basic}, $A^*$ contains a subHopf algebra of dimension 27
isomorphic to ${\mathcal B}(W) \#k [\Theta]$ and this subHopf
algebra is denoted $C(\overline{1})$. We compute the matrix
coalgebras that constitute $\gamma C(\overline{1})$ for $\gamma \not
\in \Omega .$ We work with the $k$-basis $w(\lambda, z)$ where
$\lambda \in \hat{\Gamma}$ and $z \in {\mathcal Z} = \{ 1, x,y, x^2,
y^2, xy, x^2y, xy^2, (xy)^2 \}.$

Then a $k$-basis for $C (\overline{\gamma})$ is $w(\gamma \chi^i,
z)$ for $0 \leq i \leq 2$ and $z \in \mathcal Z$ as above. As in the
proof of Theorem \ref{rank2}, since the only elements of $\mathcal
Z$ that commute with elements of $\Gamma$ are $1, xy$ and $(xy)^2$,
then any candidate for a matrix coalgebra basis element $e(i,i)$
containing $\gamma$ as a summand must be a linear combination of
$\gamma, w(\gamma, xy)$ and $w(\gamma, (xy)^2)$. Let $\eta = \gamma
(g^3 -1) \neq 0$. Supplying extra detail in this first example, we
note that
\begin{eqnarray*}
&& \triangle (\gamma) = \gamma \otimes \gamma + \eta [
w(\gamma, x)\otimes  w( \gamma \chi^2, x^2) + w(\gamma,x^2)
\otimes w (\gamma \chi, x)\\
 &+& w( \gamma, y) \otimes w (\gamma
\chi, y^2) + w(\gamma, y^2) \otimes w(\gamma \chi^2, y) ] +
\eta^2 [ q w (\gamma, xy^2) \otimes w(\gamma \chi, x^2y)\\
 &+& qw
(\gamma, x^2 y) \otimes w (\gamma \chi^2, xy^2)  + w(\gamma, xy)
\otimes w(\gamma, (xy)^2) + w(\gamma, (xy)^2) \otimes w(\gamma, xy)
]; \end{eqnarray*} \begin{eqnarray*} &&\triangle (w(\gamma, xy)) =
\gamma \otimes w (\gamma, xy) + w (\gamma, x) \otimes w(\gamma
\chi^2,y) + qw (\gamma, y) \otimes
w(\gamma \chi, x) \\
&+& \eta [w (\gamma, x^2) \otimes w(\gamma \chi, x^2y) + q^2
w(\gamma, y^2) \otimes w (\gamma \chi^2, xy^2) +   w(\gamma,
xy^2) \otimes w (\gamma \chi, y^2)\\
& +&
 q^2 w(\gamma, x^2 y)
\otimes w(\gamma \chi^2,x^2)] + \eta^2 w(\gamma, (xy)^2) \otimes
w(\gamma, (xy)^2) + w(\gamma, xy) \otimes \gamma;
\end{eqnarray*} \begin{eqnarray*} &&\triangle (w (\gamma, (xy)^2) = \gamma \otimes w (\gamma,
(xy)^2) + w (\gamma, xy) \otimes w(\gamma, xy) +
qw (\gamma,y) \otimes w(\gamma \chi, x^2y) \\
&+&   qw (\gamma,
xy^2) \otimes w (\gamma \chi, x)
+  q^2 w(\gamma, x^2) \otimes w(\gamma \chi, y^2) + q^2
w(\gamma, x) \otimes w (\gamma \chi^2, xy^2) \\
& +& w(\gamma, y^2) \otimes w(\gamma \chi^2,x^2) + q^2 w (\gamma,
x^2y) \otimes w(\gamma \chi^2, y)+ w(\gamma, (xy)^2) \otimes
\gamma.
\end{eqnarray*}

Now define $e(1,1) = \gamma + \eta^{2/3} w(\gamma, xy)
+ \eta^{4/3} w(\gamma, (xy)^2)$ as suggested by Remark \ref{choosee}.
By the computations above,
 $\triangle(e( 1,1) = e(1,1)\otimes
e( 1,1) + e(1,2) \otimes e(2,1) + e(1,3) \otimes e(3,1), $
 where
\begin{eqnarray*}
e(1,2)& =& \eta^{1/3} w(\gamma, x) + \eta^{2/3}
w(\gamma, y^2)+ q^2
\eta w(\gamma, x^2y); \\
e(2,1) &=&
\eta^{2/3} w(\gamma
\chi^2,x^2) + \eta^{1/3} w(\gamma \chi^2,y) +
q^2 \eta w (\gamma \chi^2, xy^2); \\
e(1,3) &=& \chi^2 e(2,1) \chi^2 = q\eta^{2/3} w(\gamma,x^2)+q^2
\eta^{1/3} w(\gamma, y)  +q^2
\eta w (\gamma, xy^2) ; \\
e(3,1) &=& \chi^2 e(1,2) \chi^2 = q^2  \eta^{1/3} w(\gamma \chi, x) +q
\eta^{2/3} w(\gamma \chi, y^2) + q^2 \eta w ( \gamma \chi, x^2 y).
\end{eqnarray*}
Now we compute $\triangle (e(3,1))$ and define $e( 3,3)  = \gamma
\chi + q \eta^{2/3} w(\gamma \chi, xy) + q^2 \eta^{4/3} w(\gamma
\chi, (xy)^2) = \chi^2 e(1,1) \chi^2.$ Similarly we have
\begin{eqnarray*}
 e(2,2)&=& \chi e(1,1) \chi= \gamma \chi^2 + q^2
\eta^{2/3} w(\gamma \chi^2,xy) +q
\eta^{4/3} w(\gamma \chi^2, (xy)^2);\\
e(2,3) &=& \chi e(1,2) \chi = q \eta^{1/3} w(\gamma \chi^2,x) +
q^2 \eta^{2/3}
w(\gamma \chi^2, y^2) + q^2 \eta w(\gamma \chi^2, x^2y);\\
e(3,2)& =& \chi e(2,1) \chi = q^2\eta^{2/3} w(\gamma \chi, x^2)+
q \eta^{1/3} w (\gamma \chi,y) +q^2 \eta w(\gamma
\chi, xy^2) .\\
\end{eqnarray*}
Straightforward computation then shows that these elements satisfy
the conditions to be a basis for a $3 \times 3$ matrix coalgebra and
that by the proof of Theorem \ref{rank2} or by easy computation,
this matrix coalgebra is not invariant under the inner action of
$\chi$.  Thus the other two matrix coalgebras constituting
$C(\overline{\gamma})$ may be obtained by using this action and
$C(\overline{\gamma}) = \oplus _{i=1}^3 {\mathcal M}^c(3,k)$. \qed
}\end{ex}

\vspace{.2cm}

\begin{ex}\label{order3II}
 {\rm This example is similar to Example \ref{order3I} but
$m=1$ and $n = -1$, i.e. $V =kv_1 \oplus kv_2$ with $v_1 \in
V^\chi_g$ and $v_2 \in V^\chi_{g^{-1}}$ where $\chi^3 = \epsilon$
and $g^3 \neq 1$. Otherwise $A, {\mathcal A}, x,y,q$ are as defined
in Example \ref{order3I}. Again we compute the matrix coalgebras
that constitute $C (\overline{\gamma})$ for $\gamma \not \in \Omega
= G(A^*).$ This time, we work with the $k$-basis for $A$ given by
$hz$, $h \in \Gamma$, $z \in {\mathcal Z} = \{ 1,x,y, x^2, y^2, xy,
x^2y, xy^2, x^2y^2  \}$  so $A^*$ has basis $w( \gamma \chi^k, x^i
y^j), 0 \leq k ,i,j \leq 2.$ Here the elements of ${\mathcal Z}$
that commute with grouplikes are $h, hxy^2$ or $hx^2y, h \in
\Gamma,$ so that we guess that $e(1,1)$ has the form $\gamma + aw
(\gamma, x^2y) + bw(\gamma, xy^2)$. Straightforward computation
yields that $a^3 =b^3 = (\eta q)^3 = \eta^3$. We choose $a =b = \eta
q$, and then find that
\begin{eqnarray*}
e(1,1) &=& \gamma + \eta q w(\gamma, x^2y) + \eta qw (\gamma,
xy^2);\\
e(2,2) &=& \gamma \chi^2 + \eta q^2 w(\gamma \chi^2, x^2y) + \eta
w(\gamma \chi^2, xy^2)= \chi e(1,1)\chi;\\
e(3,3) &=& \gamma \chi + \eta w (\gamma \chi, x^2y) + \eta q^2
w(\gamma \chi, xy^2) = \chi^2 e(1,1) \chi^2;\\
e(1,2) &=& \eta w(\gamma,x) + \eta q^2 w(\gamma,y) + \eta^2
w(\gamma, x^2y^2);\\
e(3,1) &=& \eta q^2 w(\gamma \chi,x) + \eta w (\gamma \chi,y) +
\eta^2 w(\gamma \chi, x^2 y^2) = \chi^2 e(1,2) \chi^2;\\
e(2,3) &=& \eta^{-1} qw(\gamma \chi^2, x) + \eta^{-1} qw (\gamma
\chi^2,y) + w(\gamma \chi^2, x^2y^2) = \eta^{-2} \chi e(1,2)
\chi;\\
e(2,1) &=& w(\gamma \chi^2 , x^2) + qw (\gamma \chi^2, y^2) + qw
(\gamma \chi^2, xy);\\
e(1,3)& =& qw(\gamma, x^2) +w(\gamma, y^2) + qw(\gamma,xy) =
\chi^2 e(2,1) \chi^2;\\
e(3,2) &=& q^2 \eta^2 w(\gamma \chi, x^2) + \eta^2 q^2 w(\gamma
\chi, y^2) + \eta^2 qw (\gamma \chi, xy) = \eta^2 \chi e(2,1)
\chi;
 \end{eqnarray*}
 form a matrix coalgebra basis for a matrix coalgebra $E = \chi E \chi =
 \chi^2 E \chi^2.$  Then $C(\overline{\gamma}) = E \oplus E \chi \oplus
 E \chi^2 \cong \oplus_{i=1}^3 {\mathcal M}^c(3,k).$ }\qed
\end{ex}

\subsection{Duals of pointed Hopf algebras of dimension 16}

Now we describe the duals of all pointed Hopf
 algebras of dimension 16 which are not group algebras or their duals.
 In fact, this means that the duals of all Hopf algebras of dimension 16
 with  coradical  a Hopf algebra are known, since the nonpointed Hopf
 algebras of dimension 16 with coradical a proper subHopf algebra were computed in
 \cite{cdm} and found to be self-dual.
 In \cite{cdr}, all pointed Hopf algebras of dimension 16 over an
 algebraically closed field $k$ of characteristic 0 were computed,
 up to isomorphism. The classification of semisimple Hopf algebras
 of dimension 16 is done in \cite{kashina}.

By \cite{cdr},  there are 43 different isomorphism classes of
pointed Hopf algebras: 14
 for the group algebras and 29 for Hopf algebras with group of
 grouplikes of order 2, 4 or 8. Of these 29 Hopf algebras, 21 are
 of the form ${\mathcal B}(V) \# k[\Gamma]$ where $V$ is a quantum
 linear space and 8  are nontrivial liftings of such ${\mathcal B}(V) \# k[\Gamma]$.
 First we list the Radford biproducts and their duals.
We use the following notation.

If $\Gamma$ is cyclic, then its generator is denoted by $c$ and
$\hat{\Gamma}$ has generator $c^*$.  If $\Gamma$ is the product of
2 cyclic groups then we denote their generators by $c$ of order
$n$
 and $d$ of order $m$ where $n \geq m$ and
$c^*,d^* $ are the elements of $\hat{\Gamma}$ which map $c$ to a
primitive $n$th root of unity and $d$ to 1, $d$ to a primitive
$m$th root of unity and $c$ to 1, respectively. If $\Gamma $ is
the product of 3 cyclic groups, we denote the generators by
$c,d,e$ and the generators of $\hat{\Gamma}$ by $c^*,d^*,e^*$
with the same conventions.

We describe the space $V= \oplus^{\dim V}_{j=1} kv_j$ with $v_j \in
V^{\chi_j}_{g_j}$ by pairs $(g_j, \chi_j)$, one pair for each $v_j$.
A list of the  $H_i = {\mathcal B} (V_i) \#k[\Gamma_i]$ and their
duals is given by the table below.
\begin{center}
\begin{tabular}{|c|l|l|l|} \hline
$H_i$& $\Gamma_i$ & $V_i$ & Dual of $H_i$ \\ \hline
$H_1$&$C_8$ & $(c,c^{*4})$ & $H_2$ \\
$H_2$&$C_8$ & $(c^4, c^*)$ & $H_1$\\
$H_3$&$C_8$ &$(c^2, c^{*2})$ & self-dual \\
$H_4$&$C_8$ & $(c^6, c^{*2}) \cong (c^2, c^{*6})$ & self-dual \\
[.10in]
$H_5$&$C_4 \times C_2$ & $(c,c^{*2})$ & $H_6$ \\
$H_6$&$C_4 \times C_2$ & $(c^2, c^*)$ & $H_5$\\
$H_7$&$C_4 \times C_2$ &$(cd, d^{*})$ & $H_8$ \\
$H_8$&$C_4 \times C_2$ & $(d,c^*d^*) \cong (c^2d, c^{*})$ & $H_7$\\
$H_9$&$C_4 \times C_2$ & $(d, d^*)$ & self-dual \\ [.10in]
$H_{10}$&$C_2 \times C_2 \times C_2$ & $(c,c^{*})$ &self-dual
\\[.10in]
$H_{11}$&$C_4$ & $(c,c^{*})$ & self-dual \\
$H_{12}$&$C_4$&$(c, c^{*3}) \cong (c^3, c^*)$ & self-dual \\
$H_{13}$&$C_4$ & $(c^2,c^*)(c^2,c^*)$ & $H_{14}$ \\
$H_{14}$&$C_4$ & $(c, c^{*2})(c, c^{*2})$ & $H_{13}$\\
$H_{15}$&$C_4$ &$(c^2, c^*)(c^2,c^{*3})$ & $H_{16}$ \\
$H_{16}$&$C_4$ & $(c, c^{*2})(c^3, c^{*2})$ & $H_{15}$\\ [.10in]
$H_{17}$&$C_2 \times C_2$ & $(c,c^*)(c,c^*)$ & self-dual \\
$H_{18}$&$C_2 \times C_2$ & $(c, c^{*})(d, d^*)$ & self-dual\\
$H_{19}$&$C_2 \times C_2$ &$(c, c^*)(c,c^{*}d^*)$ & $H_{20}$ \\
$H_{20}$&$C_2 \times C_2$ & $(c, c^*)(cd, c^{*})$ or
$(c,c^*d^*)(d, c^*d^*)$ & $H_{19}$\\ [.10in]
$H_{21}$&$C_2$ & $(c,c^*)(c,c^*)(c, c^*)$ & self-dual \\
 \hline
\end{tabular}

\end{center}

Note that for some classes, such as $H_4$ and $H_8$, we show two
choices for $g_i$ and $\chi_i$. To see that the choices for $H_8$
are isomorphic, use the automorphism $\phi$ of $\Gamma$ given by
$\phi (c) = c, \phi (d) =c^2d$. The choice $g_8 =d, \chi_8 =c^*d^*$
is natural from Corollary \ref{cordual}; the choice $g_8 =c^2d,
\chi_8 =c^*$ is mentioned because it corresponds to the listing of
the classes in \cite{cdr}. To confirm the isomorphism for $H_4$, use
the automorphism $\phi$ of $\Gamma$ given by $\phi (c) = c^3$. As
noted in \cite{cdr}, $H_{10} \cong T \otimes k[C_2 \times C_2]$
where $T$ is the 4-dimensional Sweedler Hopf algebra and so it is
obvious that $H_{10}$ is self-dual.

Now we compute  some  duals before summarizing information on the
structure of the duals of   the nontrivial liftings. Although
  Theorem \ref{rank2} does not apply here since $r$ is not
an odd integer,   the approach is similar.

 \begin{ex}\label{a142}  {\rm Let $A$ be the lifting of $H_{14}$ with   matrix
 ${\mathcal A}= I_2$, i.e. $A$ is the
Hopf algebra $H_3$ of \cite[Section 2]{beattie} or the Hopf
algebra $H_1(1)$ from \cite[Section 3]{cdr}. Thus $\Gamma= C_4 =
<c>$, $\hat{\Gamma} = < c^*>$ and $V=V^{c^{*2}}_c =kx \oplus ky$
with $xy =-yx$ and $x^2 = y^2 = c^2 -1$. We describe the dual
$A^*$.
 In this case, $\Omega = \{ 1, c^{*2}\}$. Then $A^*$ contains an
8-dimensional Hopf subalgebra $C(\overline 1)$ isomorphic to
${\mathcal B} (W) \# k [\Omega]$ where $W = kw \oplus kz =
W^c_{c^{*2}}$.  By Lemma \ref{multiplication}, as an algebra, $A^*
\cong {\mathcal B}
 (W) \# k \hat{[\Gamma]}$.  We study the 8-dimensional coalgebra
 $C(\overline{c^*}) = C(\overline{c^{*3}})$ by finding a
 matrix subcoalgebra in $  C(\overline{c^*})$
  explicitly. Note that $c^*(c^2-1) =-2$, let $i$ be a primitive 4th root of 1 and let
   $E = \oplus_{1 \leq i,j \leq 2}
 k e(i,j)$ with

 \begin{eqnarray*}
  e(1,1) &=& c^*-2i w(c^*, xy); \\
 e(1,2) &=& w(c^*,x) -
  iw(c^*,y);\\
  e(2,1) &=& -2 [w (c^{*3},x) + iw (c^{*3},y)];  \\
  e(2,2) &=&
  c^{*3}+2i
  w(c^{*3},xy).
\end{eqnarray*}
  Then checking that $\{ e(1,1),
e(1,2), e(2,1), e(2,2)\}$
 is a
matrix coalgebra basis is a straightforward computation. As noted in
Remark \ref{requal2}, we have $c^{*2}Ec^{*2} = E.$ But $E \neq F =
c^{*2}E  =  Ec^{*2} \cong {\mathcal M}^c(2,k)$. Here $F$ has matrix
coalgebra basis
 $ \{ f(i,j)| 1 \leq i,j \leq 2\} = \{ e(i,j) c^{*2} | 1 \leq i,j
 \leq 2\} = \{ c^{*3} - 2i w(c^{*3}, xy), w(c^{*3}, x) -
 iw(c^{*3},y), -2 [w(c^*,x) +iw(c^*,y)], c^* +2i w(c^*,xy) \}$.
Thus $A^*$ has coradical isomorphic to $k[C_2] \oplus _{i=1}^2
{\mathcal M}^c(2,k)$.

Now let $B$ be the lifting of $H_{16} $ with matrix ${\mathcal A} =
I_2$ and note that as algebras $A = B$, i.e. both are generated by
$C_4, x$ and $y$ with the same commutation relations.  Thus as
coalgebras, their duals are isomorphic, i.e., $B^*$ has coradical
isomorphic to $k[C_2] \oplus _{i=1}^2 {\mathcal M}^c(2,k)$ and $E$
and $F$ are the same as in $A^*.$ The algebra structure of $B^*$ is
determined in  Lemma \ref{multiplication}. }\qed
\end{ex}

It remains to describe the dual of the lifting of $H_{20}$ with
nondiagonal matrix ${\mathcal A}.$
 Note that this Hopf algebra $A$
is the Hopf algebra discussed in Example \ref{E(n)}, i.e. the
Drinfel'd double of the Sweedler 4-dimensional Hopf algebra. It is
the first lifting we have considered with a nondiagonal matrix
${\mathcal A} = \left [ \begin{array}{ll} 0& 1\\
1&0
\end{array} \right ]$.

\begin{ex}\label{a20}{\rm
Let $A$ be as above. Then $A$ is generated by $k[\Gamma]$ where
$\Gamma = <c> \times <d> \cong C_2  \times C_2,$ the
$(1,c)$-primitive $x$ and the $(1,cd)$-primitive $y$. Also $hx
=c^* (h) x h$ and $hy = c^* (h) yh$ for all $h \in \Gamma$ where
$c^* (c) =d^*(d) =-1$ and $c^*(d) = d^*(c) = +1$. In $A$, $x^2
=y^2=0$, and $xy +yx =d-1$.

Then $\Omega = \{ 1,c^*\}$ is  the group of grouplikes  of $A^*$ and
$A^*$ contains a subHopf algebra $C(\overline{1})$ of dimension 8
isomorphic to ${\mathcal B}(W) \# k [\Omega]$ where $W = kw_1 \oplus
kw_2$ with $w_1 = w(c^*,x)$ and $w_2 = w(c^*,y)$. Then $w_1 \in
W^c_{c^*}$, $w_2 \in W^{cd}_{c^*}$ and $w_1^2 = w_2^2 = 0.$

Note that the standard $k$-basis we have been using for $A$ is
$\{ g, gx, gy, gxy | g \in \Gamma\}$ but the set $\{ g, g(x+y),
g(x-y), g(x+y)(x-y) | g\in \Gamma \}$ is also a basis. Let $x'
=x+y, y' = x-y$. Then the set $\{ w(\gamma, 1) =\gamma, w(\gamma,
x'), w(\gamma, y'), w(\gamma, x'y') | \gamma \in \hat{\Gamma} \}$
is a $k$-basis for $A^*$ where $ w(\gamma,z)(hz) = \gamma (h)$ and
$w(\gamma, z)$ is zero on all other basis elements as usual. Note
that $2w(\gamma, x') = w(\gamma, x) +w(\gamma,y)$ and also that
$w(c^*, x'), w(c^*,y')$ are $(\epsilon, c^*)$-primitive. Thus
$C(\overline{d}^*)$ has $k$-basis $\{ w(d^* \gamma, z)| \gamma
\in \{ 1,c^*\}, z \in \{ 1, x',y', x'y'\}\}.$ Note that $x'^2= -
y'^2 =xy +yx =d-1$. Also $x'y' =yx -xy = - y'x'$. Then $(x'y')^2
=x' (-x'y')y' = (d-1)^2$. The set
\begin{eqnarray*}
e(1,1) &=& d^* + 2w(d^*, x'y'); \\
e(1,2) &=& \sqrt{2} w(d^*, x') + \sqrt{2} w(d^*, y') \\
e(2,1) &=& - \sqrt{2} w(c^*d^* ,x') + \sqrt{2} w(c^*d^*, y') ; \\
e(2,2) &=& c^*d^* - 2w(c^*d^*, x'y')
\end{eqnarray*}
is a matrix coalgebra basis for a matrix coalgebra $E$. A sample
computation is:
\begin{eqnarray*}
\triangle e(1,1) &=& \triangle d^* + 2 \triangle w(d^*, x'y')\\
&=& [d^* \otimes d^* + (-2) w(d^*,x') \otimes w(c^*d^*, x')
 + (2)w(d^*,y') \otimes w(c^*d^* ,y') \\
 &+& 4w(d^*, x'y') \otimes w
(d^*,x'y')]
 + 2 [d^* \otimes w (d^*, x'y') + w(d^*, x'y') \otimes d^*\\
& +&w(d^*, x') \otimes w(c^*d^*,y') -w (d^*,y') \otimes w(c^*d^*. x')]\\
&=& d^* \otimes [d^* + 2w (d^* , x'y')]
+ 2w (d^*, x'y') \otimes [d^* + 2w (d^*, x'y')] \\
&+& \sqrt{2} w
(d^*, x') \otimes [ -\sqrt{2} w(c^*d^*, x') + \sqrt{2} w(c^*d^*, y')]\\
&+& \sqrt{2} w(d^*, y') \otimes [ \sqrt{2} w(c^*d^*, y') -
\sqrt{2} w
(c^*d^*, x')]\\
&=& e(1,1) \otimes e(1,1) + e(1,2) \otimes e(2,1).
 \end{eqnarray*}

Note that the first multiplication formula in Lemma
\ref{multiplication}(ii) does not apply to products such as $c^*
* w(d^*, x')$ because $x'$ is not the product of skew-primitives.
However, the formulas do apply to products such as
$c^**w(d^*,x'y')$ since $\triangle(x'y') = \triangle (yx - xy)$
contains a term $d \otimes x'y'.$ Thus we see that $c^* e(1,1)c^*
= e(1,1)$ and so $c^*Ec^* = E.$  Or we could recall the fact that
$c^*Ec^* = E$ from Remark \ref{requal2}.

\par Also $c^*E = Ec^* = F \neq E$ is a matrix coalgebra  with matrix
coalgebra basis $f(i,j) = e(i,j)c^*$.

\par Thus we have shown that as a coalgebra $A^*  \cong
C(\overline{1}) \oplus C(\overline{d^*})$ where
 $C(\overline{d^*})$ is the direct sum of
two 4-dimensional matrix coalgebras, so that corad $(A^* ) \cong k
[C_2] \oplus C (\overline{d}^*)$. As an algebra, $A^* $ is generated
by $\hat{\Gamma}, w_1 = w(c^*, x)$ and $w_2 =w(c^*,y)$ where for
$\gamma \in \hat{\Gamma}$, $\gamma w_1 = \gamma (c) w_1 \gamma
\mbox{ and } \gamma w_2 = \gamma (cd) w_2 \gamma; \quad w_1 w_2 =-
w_2 w_1;\quad w^2_1 = w^2_2 =0$. Alternatively, as an algebra $A^*$
is generated by $\hat{\Gamma}$ and $w'_1 = w(c^*, x'), w'_2 = w(c^*,
y')$ but $\gamma (kw'_i) \gamma^{-1} \neq kw'_i$.} \qed
\end{ex}

Now we list the duals of the 8 classes of pointed Hopf algebras of
dimension 16 which are nontrivial liftings. In the table below,
$A_i$ is the nontrivial lifting of the Hopf algebra $H_i$ with
matrix ${\mathcal A}_i$,  where  $D=\left [
\begin{array}{ll} 1& 0\\ 0&0 \end{array} \right ]$ and $E=\left [
\begin{array}{ll} 0& 1\\ 1&0 \end{array} \right ] $.

\begin{center}

\begin{tabular}{|c|l|l|l|l|l|} \hline
$A_i$ & $G(A_i)$ & $V_i$ & ${\mathcal A}_i$ & $A_i^*$ & $corad(A_i^*)$\\
\hline
&&& &&\\
 $A_1$ & $C_8$ & $(c, c^{*4})$ & $I_1$ & Th. \ref{rank1}&
$k[C_2] \oplus \oplus^3_{i=1} {\mathcal M}^c (2,k)$\\[.10in]

$A_5$ & $C_4 \times C_2  $ &$(c, c^{*2})$ & $I_1$ &Th. \ref{rank1}
& $k[C_2 \times C_2] \oplus \oplus^2_{i=1} {\mathcal M}^c (2,k)$\\[.10in]

$A_7$ & $C_4 \times C_2 $ & $(cd,d^*)$ & $I_1$ & Th. \ref{rank1}
& $k[C_2 \times C_2] \oplus \oplus^2_{i=1} {\mathcal M}^c (2,k)$\\[.10in]

$A_{14,1}$ & $C_4 $ & $(c,c^{*2}), (c, c^{*2})$ & $D$&
Cor.\ref{quotient} &$ k[C_2] \oplus  {\mathcal M}^c (2,k)$\\[.10in]

 $A_{14,2}$ & $C_4  $ &
$(c,c^{*2}), (c, c^{*2})$ & $I_2$ & Ex. \ref{a142}
&$ k[C_2] \oplus \oplus^2_{i=1} {\mathcal M}^c (2,k)$ \\[.10in]

$A_{16,1}$ & $C_4  $ & $(c,c^{*2}), (c^3, c^{*2})$ & $ D
  $& Cor.\ref{quotient}
&$ k[C_2] \oplus {\mathcal M}^c (2,k)$\\[.10in]
$A_{16,2}$ & $C_4  $ & $(c,c^{*2}), (c^3, c^{*2})$ & $I_2 $ & Ex.
\ref{a142}
&$ k[C_2] \oplus \oplus^2_{i=1} {\mathcal M}^c (2,k)$\\[.10in]
$A_{20}$ & $C_2  \times C_2  $ & $(c,c^*), (cd, c^{*})$ & $ E
  $& Ex. \ref{a20}
& $k[C_2] \oplus \oplus^2_{i=1} {\mathcal M}^c (2,k)$ \\
&&&&&\\
\hline \end{tabular}
\end{center}

\subsection{Duals of an infinite family of dimension 32}

 Liftings of quantum linear spaces have provided many
 counterexamples to Kaplansky's 10th conjecture (see \cite{asp3},
 \cite{bdg}, \cite{gelaki}). The examples of smallest dimension
 are families of Hopf algebras of dimension 32 with coradical
 $k[\Gamma]$ where $\Gamma$ is an abelian group of order 8
 (see \cite{beattie}, \cite{grana32}). We compute the duals for the
 family in \cite{beattie}. Let $\Gamma = <c>
 \cong C_8$ and $\hat{\Gamma} = <c^*>$ as usual with $c^*(c) =q$,
 a primitive 8th root of 1. Let $V= kv_1 \oplus kv_2$ with $v_1
 \in V_c^{c^{*4}}, v_2 \in V_{c^5}^{c^{*4}}$ and let $A(\lambda)$ be
 the lifting of ${\mathcal B}(V) \#k[\Gamma]$ with matrix $\left [
 \begin{array}{cc} 1 & \lambda\\ \lambda &1 \end{array} \right ]$.
 Then $A(\lambda) \not \cong A(\lambda ')$ unless
 $\lambda ' = \pm \lambda$ by \cite{beattie},
 so that the $A(\lambda)$ give an
 infinite family of nonisomorphic  Hopf algebras of dimension 32.
 Since the original submission of this note, a paper by Etingof
 and Gelaki has appeared where it is shown that this family contains
 infinitely many twist equivalence classes
   \cite[Corollary 4.3]{eg}
 ,
 providing a counterexample to a conjecture of Masuoka
 \cite{masuoka}. Actually Etingof and Gelaki are studying  triangular Hopf algebras
 $A(G,V,u,B)$ in the notation of \cite{eg} (see also
  \cite{aeg} and \cite{eg1}) and the $A(\lambda)$ are duals of such Hopf algebras.
 Here, we compute the coradical of the dual of $A(\lambda)$ using the
 same type of straightforward computations as in the previous
 examples.

 \begin{ex} \label{dim32} {\rm
For $A = A(\lambda)$ as above, let $x,y$ be the liftings of $v_1$
and $v_2$. Then $x^2 =y^2 =c^2-1$ and $xy + yx = \lambda (c^6 -1)$.
Then $\Omega =G(A^*) =\{ \epsilon, c^{*4}\},$ and $A^*$ contains a
subHopf algebra $C(\overline{1})$ of dimension 8 isomorphic to
${\mathcal B}(W) \# k [\Omega]$ where $W =kw_1 \oplus kw_2$, $w_1
=w(c^{*4}, x)$ and $w_2 = w(c^{*4},y)$. We  determine the coalgebra
structure of $C(\overline{\gamma}), \gamma \not \in \Omega$.

As in Example \ref{a20}, we use a different $k$-basis for $A$. Let
$x' =x+y, y'= x-y$ and then $A$ has $k$-basis $hz$ where $h \in
\Gamma, z \in {\mathcal Z} =\{ 1, x', y', x'y' \}$. Then $x'^{2} =
2( c^2-1) + \lambda (c^6 -1), y'^{2} =2(c^2-1) - \lambda (c^6 -1)$
and $x'y' = yx -xy =- y'x'$. We use the corresponding basis $\{
w(\gamma, z) | \gamma \in \hat{\Gamma}, z \in {\mathcal Z}\}$ for
$A^*$. Let $\gamma \not \in \Omega$ and denote $\alpha = \gamma
(2(c^2-1) + \lambda (c^6-1))$, $\beta = \gamma (2(c^2-1) -\lambda
(c^6 -1)).$

Case i). Suppose first that $\alpha \neq 0$ and $\beta \neq 0.$ Then
$E = \oplus_{1 \leq i,j \leq 2} ke(i,j) \cong {\mathcal M}^c (2,k)$
where
\begin{eqnarray*}
e(1,1) &=& \gamma + i \sqrt{\alpha} \sqrt{\beta} w(\gamma, x'y');\\
e(1,2) &=& - \sqrt{\alpha} w(\gamma, x') + i \sqrt{\beta}
w(\gamma, y');\\
e(2,1) &=& - \sqrt{\alpha} w(\gamma c^{*4}, x') -i \sqrt{\beta}
w(\gamma c^{*4}, y');\\
e(2,2) &=& \gamma c^{*4} -i \sqrt{\alpha} \sqrt{\beta} w(\gamma
c^{*4}, x'y').
\end{eqnarray*}

(Here $i^2 =-1$.) Note that $E \neq E c^{*4}$ so $E \neq F =Ec^{*4}
\cong {\mathcal M}^c (2,k)$ also and $C(\overline{\gamma}) \cong
{\mathcal M}^c (2,k) \oplus {\mathcal M}^c (2,k)$.

Case ii). One of $\alpha, \beta$ is 0 and the other is nonzero.
Suppose that $\alpha \neq 0$ but $\beta =0$. Then let $E =
\oplus_{1 \leq i,j \leq 2} ke(i,j)$ where
\begin{eqnarray*}
 e(1,1) = \gamma    &; &  e(1,2) =- \sqrt{\alpha} w(\gamma, x');\\
e(2,1) = - \sqrt{\alpha} w(\gamma c^{*4}, x')&;&  e(2,2) = \gamma
c^{*4};
\end{eqnarray*}
 and it is easy to check that $E \cong {\mathcal M}^c (2,k)$.
Here $E=c^{*4} E = Ec^{*4}.$ By the same argument as in the proof
of Theorem \ref{rank2}, any matrix coalgebra basis element
$e(i,i)$ must be of the form $\gamma + \zeta w(\gamma, x'y')$.
 But $\zeta 'w (\gamma, x'y') \otimes w(\gamma, x'y')$ does not
occur in $\triangle w$
 for any basis element $w$ and so $C(\overline{\gamma})$
contains only one matrix coalgebra. It is easy to see that
$C(\overline{\gamma}) $ has coradical filtration of length 1,
i.e. $C(\overline{\gamma}) = C(\overline{\gamma})_1$. The argument
for the case $\alpha = 0, \beta \not = 0$ is the same.

Case iii). The case $\alpha = \beta =0$ cannot occur since then
$\gamma(c)^2 = 1$ so that $\gamma \in \{ \epsilon, c^{*4} \} =
\Omega$.

Now suppose $\gamma =c^*$. Then $\alpha =c^* (c^2-1) c^* (2 +
\lambda (c^4 +c^2+1)) = (q^2-1) (2+q^2 \lambda)$ since $q^4 =-1$. So
$\alpha =0$ if and only if $2+q^2 \lambda =0$. Similarly $\beta =0$
if and only if $2-q^2 \lambda =0$. If $\gamma =c^{*3}$ then $\alpha
=0$ if and only if $2+q^6 \lambda =0$ and $\beta=0$ if and only if
$2-q^6 \lambda =0$. Thus if $q^2 \lambda = \pm 2$, then
$C(\overline{c^{* }})$ and $C(\overline{c^{*3}})$ each have
coradical isomorphic to ${\mathcal M}^c (2,k).$ Otherwise both
$C(\overline{c^{* }})$ and $C(\overline{c^{*3}})$ are cosemisimple,
isomorphic to ${\mathcal M}^c (2,k) \oplus {\mathcal M}^c (2,k)$.

A similar argument shows that for $\gamma = c^{*2}$, $\alpha =0$ if
and only if $\lambda =-2$ and $\beta =0$ if and only if $ \lambda
=2$. Thus $C(\overline{c^{*2}})$ is cosemisimple unless $\lambda =
\pm 2$ and then it has coradical isomorphic to ${\mathcal M}^c
(2,k)$.

We summarize in the following table:
\begin{center}
\begin{tabular}{|c|l|} \hline
$\lambda$ & $corad (A(\lambda)^*)$ \\  \hline \phantom & \\ $\pm
2q^2$ & $k
[C_2] \oplus \oplus^4_{i=1} {\mathcal M}^c (2,k)$ \\
$\pm 2$ &$k
[C_2] \oplus \oplus^5_{i=1} {\mathcal M}^c (2,k)$ \\
$k - \{ \pm 2q^2, \pm 2 \}$ & $k [C_2] \oplus \oplus^6_{i=1}
{\mathcal M}^c (2,k)$ \\ \hline \end{tabular}
\end{center}

Note that in the last row, $\lambda$ may be 0, and then the
result follows from a simpler argument with $\alpha = \beta$. \qed
 }
 \end{ex}

{\bf Acknowledgement:} Many thanks to N. Andruskiewitsch for
pointing out more general statements for Theorem 2.2 and Theorem
2.5, and for his many insightful comments, especially  on the
material in Section 2. Also thanks to S. D\u{a}sc\u{a}lescu and S.
Westreich for comments on the original version of Theorem 2.2.


\begin{thebibliography}{99}

\bibitem{ad} N. Andruskiewitsch and S. D\u{a}sc\u{a}lescu, On finite quantum groups
at -1,  Algebr. Represent. Theory  8  (2005),  no. 1, 11--34.

\bibitem{aeg} N. Andruskiewitsch, P. Etingof, S. Gelaki,
Triangular Hopf algebras with the Chevalley property.  Michigan Math. J.  49  (2001),  no. 2, 277--298.

\bibitem{ag} N. Andruskiewitsch and M. Gra\~{n}a, Braided Hopf algebras over non-abelian finite groups,
    Bolet\'{i}n de la Acad. Nac. Cs. C\'{o}rdoba, Argentina 63
    (1999), 45-78.


\bibitem{assurvey}N. Andruskiewitsch and H.-J. Schneider,
  Pointed Hopf algebras.  New directions in
Hopf algebras,  1--68, Math. Sci. Res. Inst. Publ., 43, Cambridge
Univ. Press, Cambridge, 2002.


\bibitem{asp3}  N. Andruskiewitsch and H.-J. Schneider,
Lifting of quantum linear spaces and pointed Hopf algebras of
order $p^3$, J. Algebra 209 (1998), 659-691.

\bibitem{asa2} N. Andruskiewitsch and H.-J. Schneider, Lifting of Nichols
algebras of type $A_2$ and pointed Hopf algebras of order $p^4$,
in ``Hopf algebras and quantum groups'', eds. S. Caenepeel and F.
van Oystaeyen, M. Dekker, 2000, pp. 1-16.

\bibitem{asadvances} N. Andruskiewitsch and H.-J. Schneider,
Finite quantum groups and Cartan matrices, Adv. Math. 154
(2000),  1--45.

\bibitem{beattie}M. Beattie, An isomorphism theorem for Ore extension Hopf
algebras, Comm. Algebra 28(2)(2000), 569-584.

\bibitem{bdg} M. Beattie, S. D\u{a}sc\u{a}lescu, L. Gr\"{u}nenfelder,
Constructing pointed Hopf algebras by Ore extensions, J. Algebra
225, (2000), 743-770.

\bibitem{bdr} M. Beattie, S. D\u{a}sc\u{a}lescu and S. Raianu,
Lifting of Nichols algebras of type $B\sb 2$. With an appendix by
the authors and Ian Rutherford.  Israel J. Math.  132  (2002),
1--28.

\bibitem{cdr} S. Caenepeel, S. D\u{a}sc\u{a}lescu, and S. Raianu,
Classifying pointed Hopf algebras of dimension 16, Comm. Algebra 28(2)(2000),
541-568.

\bibitem{cdm}C\u alinescu, C.; D\u asc\u alescu, S.; Masuoka, A.; Menini, C.
Quantum lines over non-cocommutative cosemisimple Hopf algebras.
J. Algebra  273  (2004),  no. 2, 753--779.

\bibitem{cm} W. Chin and I. Musson, The coradical filtration for quantized
enveloping algebras, J. London Math. Soc. 53
(1996), 50-67.

\bibitem{drin} V. Drinfel'd, Quantum
groups, Proceedings of the International Congress of
Mathematicians, Vol. 1, 2 (Berkeley, Calif., 1986), 798--820,
Amer. Math. Soc., Providence, RI, 1987.

\bibitem {eg}
P. Etingof and S. Gelaki,  On families of triangular Hopf
algebras, IMRN \textbf{2002:14} (2002), pp.  757-768.

\bibitem {eg1}P. Etingof and S. Gelaki, On cotriangular Hopf
algebras, Amer. J. Math. 123 (2001) , 699-713.

\bibitem{gelaki} S. Gelaki, Pointed Hopf algebras and
Kaplansky's 10th conjecture, J. Algebra 209 (1998), 635-657.


\bibitem{grana} M. Gra\~{n}a, A freeness theorem for Nichols algebra, J. Algebra,
231 (2000), 235-257.

\bibitem{grana32} M. Gra\~{n}a, Pointed Hopf algebras of dimension 32,
Comm. Algebra 28 (2000), no. 6, 2935--2976.
\bibitem{kashina} Y. Kashina, Classification of semisimple Hopf
algebras of dimension 16, J. Algebra 232 (2000), 617-663.

\bibitem{kauffrad} L.H.Kauffman and D.E.Radford, A necessary and
sufficient condition for a finite dimensional Drinfel'd double to
be a ribbon Hopf algebra, J. Algebra 159 (1993), 98-114.

\bibitem{lam} T.Y. Lam, A First Course in Noncommutative Rings, Springer GTM 131.

\bibitem{lusztig} G. Lusztig, Introduction to quantum groups,
Birkh\"auser, 1993.

\bibitem{majid} S. Majid,
  Foundations of Quantum Group Theory.  Cambridge Univ. Press, 1995.

\bibitem{masuoka} A. Masuoka, Defending the negated Kaplansky
conjecture, Proc.   Amer. Math. Soc. 129 (2001), 3185-3192.

\bibitem{mont} S. Montgomery, Hopf algebras and their actions on rings,
AMS (1993), CBMS {\bf 82}.

\bibitem{adrianacomm} A. Nenciu,
 Quasitriangular structures for a class of pointed Hopf algebras
constructed by Ore extensions, Comm. Algebra 29 (2001), 1959-1981.

\bibitem{adrianajalg} A. Nenciu,
Quasitriangular pointed Hopf algebras constructed by Ore
extensions, preprint.

\bibitem{nichols} W.D. Nichols, Bialgebras of type 1,
Comm. Algebra {\bf 6} (1978), 1521-1552.


\bibitem{pvo} F. Panaite and F. Van Oystaeyen, Quasitriangular
structures for some pointed Hopf algebras of dimension $2^n$,
Comm. Algebra 27(10)(1999), 4929-4942.


\bibitem{rad} D. E. Radford, On the coradical of a finite-dimensional Hopf algebra,
Proc. Amer. Math. Soc. 53 (1975), 8-15.

\bibitem{radminqt} D. E. Radford, Minimal quasitriangular Hopf algebras,
J. Algebra 157(1993), 285-315.

\bibitem{radkauff} D.E. Radford, On Kauffman's knot invariants
arising from finite dimensional Hopf algebras, in ``Advances in
Hopf algebras'', pp205-266, J. Bergen and S. Montgomery (eds),
Lecture Notes in Pure and Appl. Math 158, Marcel Dekker, New
York, 1994.



\bibitem{Sweedler} M.E. Sweedler, Hopf Algebras, Benjamin, New York, 1969.














\end{thebibliography}
\end{document}